\documentclass[11pt]{article}
\usepackage[margin=1.0truein]{geometry}
\usepackage{setspace}
\usepackage{lscape}
\usepackage{pdflscape}

\usepackage{xcolor}

\usepackage{amsmath,amssymb,amsthm} % AMS系，ほぼ必須
\usepackage{wrapfig} % wrapfigure 用

\usepackage{mathrsfs} % 花文字フォント (\mathscr)
\usepackage{url} % url 用フォント (\url)
\usepackage{latexsym} % マニアックな数学記号用
\usepackage{bm} % 斜体太字フォント (\bm)

\usepackage{mathtools}
\usepackage{booktabs,multirow}
\usepackage{array}
\newcolumntype{P}[1]{>{\centering\arraybackslash}p{#1}}
\usepackage{diagbox}
\usepackage{threeparttable}
\usepackage{siunitx} % for scientific notation
\usepackage{algorithm}
\usepackage[noend]{algpseudocode}

\algnewcommand{\Break}{\textbf{break}}
\algnewcommand{\algorithmicgoto}{\textbf{go to} Line}
\algnewcommand{\Goto}[1]{\algorithmicgoto~\ref{#1}}
\algblockdefx{MRepeat}{EndRepeat}{\textbf{repeat}}{}
\algnotext{EndRepeat}
\usepackage{natbib}
\usepackage{thmtools,thm-restate}

\usepackage{hyperref}
\usepackage{cleveref}
\expandafter\def\csname ver@etex.sty\endcsname{3000/12/31} % autonum のバグ対策

\usepackage{breakurl}

\usepackage{multicol}
\usepackage{wasysym} % 顔文字
\usepackage{adjustbox}
\usepackage{xparse} % 複雑なマクロ定義用
\usepackage{pbox} % pbox コマンド用
\usepackage{authblk} % 著者名
\usepackage{xspace} % newcommand マクロ直後に適切にスペースを入れる

\usepackage[subrefformat=parens]{subcaption}
\usepackage{enumitem} % enumerate package より詳細な設定
\hypersetup{
  colorlinks,
  linkcolor={green!35!black},
  citecolor={green!35!black},
  urlcolor={blue!50!black}
}

\usepackage{times}
\usepackage{soul}
\usepackage{secdot}
\usepackage{cases}
\usepackage{subcaption}
\usepackage{caption}

\declaretheoremstyle[
  shaded={bgcolor=gray!15},
]{thmsty}

\declaretheorem[
  name=Theorem,
  refname={Theorem,Theorems},
  style=thmsty,
]{theorem}

\declaretheorem[
  name=Proposition,
  refname={Proposition,Propositions},
  style=thmsty,
  sibling=theorem,
]{proposition}
\declaretheorem[
  name=Lemma,
  refname={Lemma,Lemmas},
  style=thmsty,
  sibling=theorem,
]{lemma}
\declaretheorem[
  name=Definition,
  refname={Definition,Definitions},
  style=thmsty,
]{definition}
\declaretheorem[
  name=Assumption,
  refname={Assumption,Assumptions},
  style=thmsty,
]{assumption}

\declaretheorem[
  name=Example,
  refname={Example,Examples},
  style=thmsty,
]{example}

\crefname{algorithm}{Algorithm}{Algorithms}
\crefname{line}{Line}{Lines}
\crefname{section}{Section}{Sections}
\crefname{appendix}{Appendix}{Appendices}
\crefname{table}{Table}{Tables}
\crefname{figure}{Figure}{Figures}
\crefname{equation}{}{}
\Crefname{equation}{Eq.}{Eqs.}

\setlist[itemize]{
  topsep=0.4\baselineskip,
  itemsep=0\baselineskip,
  leftmargin=1.5em,
}

\setlist[enumerate]{
  font=\upshape,
  label=(\alph*),
  ref=(\alph*),
  topsep=0.4\baselineskip,
  itemsep=0\baselineskip,
  leftmargin=2em,
}

\newlist{enuminasm}{enumerate}{1} % asm 環境内で enumerate を使うときの環境
\setlist[enuminasm]{
  font=\upshape,
  label=(\alph*),
  ref=\theassumption(\alph*),
  topsep=0.4\baselineskip,
  itemsep=0\baselineskip,
  leftmargin=2em,
}
\crefalias{enuminasmi}{assumption}

\newlist{enuminthm}{enumerate}{1}
\setlist[enuminthm]{
  font=\upshape,
  label=(\alph*),
  ref=\thetheorem(\alph*),
  topsep=0.4\baselineskip,
  itemsep=0\baselineskip,
  leftmargin=2em,
}
\crefalias{enuminthmi}{theorem}

\newlist{enuminlem}{enumerate}{1}
\setlist[enuminlem]{
  font=\upshape,
  label=(\alph*),
  ref=\thelemma(\alph*),
  topsep=0.4\baselineskip,
  itemsep=0\baselineskip,
  leftmargin=2em,
}
\crefalias{enuminlemi}{lemma}

\newlist{enuminprop}{enumerate}{1}
\setlist[enuminprop]{
  font=\upshape,
  label=(\alph*),
  ref=\theproposition(\alph*),
  topsep=0.4\baselineskip,
  itemsep=0\baselineskip,
  leftmargin=2em,
}
\crefalias{enuminpropi}{proposition}

\newlist{enumincond}{enumerate}{1}
\setlist[enumincond]{
  font=\upshape,
  label=(\alph*),
  ref=\thecondition(\alph*),
  topsep=0.4\baselineskip,
  itemsep=0\baselineskip,
  leftmargin=2em,
}
\crefalias{enumincondi}{condition}

\DeclareMathOperator*{\argmax}{argmax}

\DeclareMathOperator*{\tsum}{\textstyle\sum}

\DeclareMathOperator{\EE}{\mathbb{E}}

\DeclarePairedDelimiterX\inner[2]{\langle}{\rangle}{{#1},{#2}}

\DeclarePairedDelimiter\bra{[}{]}
\DeclarePairedDelimiterX\Set[2]{\{}{\}}{\mspace{2mu}{#1}\;\delimsize|\;{#2}\mspace{2mu}}
\DeclarePairedDelimiterX\Prn[2]{(}{)}{\mspace{2mu}{#1}\;\delimsize|\;{#2}\mspace{2mu}}
\DeclarePairedDelimiterX\Bra[2]{[}{]}{\mspace{2mu}{#1}\;\delimsize|\;{#2}\mspace{2mu}}

\newcommand{\R}{\mathbb R}

\newcommand{\bx}{\bm{x}}

\newcommand{\ba}{\bm{a}}
\newcommand{\bp}{\bm{p}}

\newcommand{\bv}{\bm{v}}

\newcommand{\bme}{\bm{e}}
\newcommand{\by}{\bm{y}}
\newcommand{\bmxi}{\bm{\xi}}

\newcommand{\bg}{\bm{g}}

\newcommand{\blambda}{\bm{\lambda}}
\newcommand{\mR}{\mathbb{R}}
\newcommand{\bmu}{\bm{\mu}}

\newcommand{\E}{\mathbb{E}}
\newcommand{\mpr}{\mathrm{Pr}}

\renewcommand{\epsilon}{\varepsilon}
\NewDocumentCommand{\exsub}{s m O{} m}{%
  \IfBooleanT{#1}{\EE_{#2}\nolimits\bra*{#4}}%
  \IfBooleanF{#1}{\EE_{#2}\nolimits\bra[#3]{#4}}%
}
\NewDocumentCommand{\ex}{s O{} m}{%
  \IfBooleanT{#1}{\EE\nolimits\bra*{#3}}%
  \IfBooleanF{#1}{\EE\nolimits\bra[#2]{#3}}%
}
\NewDocumentCommand{\cex}{s O{} m m}{%
  \IfBooleanT{#1}{\EE\nolimits\Bra*{#3}{#4}}%
  \IfBooleanF{#1}{\EE\nolimits\Bra[#2]{#3}{#4}}%
}
\usepackage{CJKutf8}

\newcommand{\email}[1]{\href{mailto:#1}{\nolinkurl{#1}}}

\usepackage{amsmath}

\usepackage{datetime}
\date{\vspace{-2.5\baselineskip}}

\author[1]{Yuya Hikima\footnote{Corresponding author. E-mail: \url{yuya-hikima@g.ecc.u-tokyo.ac.jp}}}
\author[1,2]{Akiko Takeda}

\affil[1]{Graduate School of Information Science and Technology, University of Tokyo, Tokyo, Japan}
\affil[2]{Center for Advanced Intelligence Project, RIKEN, Tokyo, Japan}

\title{Stochastic Approach for Price Optimization Problems with Decision-dependent Uncertainty}

\begin{document}
\maketitle

\begin{abstract}
Price determination is a central research topic of revenue management in marketing.
The important aspect in pricing is controlling the stochastic behavior of demand, and the previous studies have tackled price optimization problems with uncertainties.
However, many of those studies assumed that uncertainties are independent of decision variables (i.e., prices) and did not consider situations where demand uncertainty depends on price.
Although some price optimization studies have dealt with decision-dependent uncertainty, they make application-specific assumptions in order to obtain an optimal solution or an approximation solution.
To handle a wider range of applications with decision-dependent uncertainty, we propose a general non-convex stochastic optimization formulation. 
This approach aims to maximize the expectation of a revenue function with respect to a random variable representing demand under a decision-dependent distribution.
We derived an unbiased stochastic gradient estimator by using a well-tuned variance reduction parameter and used it for a projected stochastic gradient descent method to find a stationary point of our problem.
We conducted synthetic experiments and simulation experiments with real data on a retail service application.
The results show that the proposed method outputs solutions with higher total revenues than baselines.
\end{abstract}

\section{Introduction}\label{sec:introduction}
Price determination is a central research topic of revenue management in marketing, and many pricing studies have targeted applications in agricultural \citep{wang2019integrating}, online retail \citep{ferreira2016analytics}, electrical power \citep{dong2017electricity}, and hospitality industries \citep{koushik2012retail}.

An important aspect in pricing is controlling the stochastic behavior of demand. 
This is because stochastic over/under demand causes a loss in many cases; for example, in road pricing, overuse of a certain road causes congestion or traffic accidents; 
in an electricity market, if demand is much lower than the available electricity supply, capital investment costs cannot be recovered.

To obtain greater profits under demand uncertainty, many of the previous studies have tackled price optimization problems with decision-independent random variables.
For example, \cite{yong2009coordinating} and  \cite{dong2017electricity} define the demand for a product/service as $d(x)+\xi$, where $x$ is price and $\xi$ is a decision-independent random variable.
\cite{correa2017posted} and \cite{chawla2010multi} assume multi-agent systems where each buyer $i$ has a random variable $v_i$ as their value for a product and purchases it when the price is below the value.
However, in practical applications, it is natural for the distribution of stochastic demand to vary with price: 
when the price of a product is close to (far from) those of competing products, it is difficult (easy) to predict the demand and its uncertainty is large (small).
Furthermore, the settings of these studies with decision-independent random variables need to use discontinuous functions to represent buyers' discrete actions (e.g., buy or leave), which makes the optimization problem difficult to solve (see Section \ref{subsec:related_independent_formlulation}).

Although some pricing studies have dealt with decision-dependent uncertainty, they assume specific demand distributions and problem settings in order to obtain an optimal solution or an approximation solution.
For example, \cite{Bertsimas2005dynamic} determine prices for multiple products produced with limited resources.
They consider demand of item $i$ at time $t$, i.e.,  $\alpha^t_i(x^t_i)+\beta^t_i(x^t_i)\xi^t_i$, where $\alpha^t_i(\cdot)$ and $\beta^t_i(\cdot)$ are given functions, $x^t_i$ is price, and $\xi^t_i$ is a random variable following a given distribution.
\cite{schulte2020price} optimize prices over multiple periods to sell a single product with a fixed unit cost $c \ge 0$, where the item's demand follows a Poisson distribution with a given intensity function $\lambda(x)$.
While these studies can find optimal or approximation solutions for their problems,
they appear to be difficult to apply to a wide range of probability distributions or problem settings (e.g., a nonlinear cost setting for selling products) due to their specific assumptions.

To resolve these issues, for general price optimization, we propose a non-convex stochastic optimization formulation that maximizes the expectation of a revenue function with respect to a random variable representing demand under a decision-dependent distribution.
Our formulation assumes that (i) the objective function is differentiable and Lipschitz continuous, (ii) the given probability density function of the random variables is differentiable and its gradient, normalized by the value of the probability density function, is bounded, and (iii) the feasible region is compact and convex.
These assumptions may seem strong, but they often hold in the price optimization literature.
Indeed, we show three application examples satisfying our assumptions (see Section \ref{sec:applications}).

The formulated problem for practical applications is generally non-convex and the dimension of the decision variables may be large.
We derive an unbiased stochastic gradient estimator of the objective function by using information on the probability density function and incorporate the estimator in a projected stochastic gradient descent method to find a stationary point of our problem.
When deriving a gradient estimator, it is important to design it so that its variance is small for fast convergence of the algorithm.
Our unbiased stochastic gradient includes a variance reduction parameter, which is inspired by \textit{baseline} technique \citep{williams1992simple,sutton2018reinforcement} in the reinforcement-learning literature. 
After confirming that the variance of the proposed stochastic gradient is bounded, we present a method for calculating the variance reduction parameter.
Then, we develop a projected stochastic gradient descent method, which converges to a stationary point by incorporating the proposed stochastic gradient and method for calculating the variance reduction parameter into a recent gradient descent algorithm \citep{ghadimi2016accelerated}. 
Moreover, we show a way of speeding up the computation of the minibatch gradient under additional assumptions that hold in applications where multiple agents make purchase decisions.

While some of the previous methods might seem applicable to our formulation, they are not suitable for the following reasons: the retraining method \citep{perdomo2020performative,mendler2020stochastic} requires strong convexity of the objective function;
the Bayesian optimization \citep{brochu2010tutorial,frazier2018tutorial} and gradient-free methods \citep{spall2005introduction,Flaxman2004online} require a huge number of evaluations of objective values, which makes it difficult to find good solutions for large-scale problems in a reasonable time.

We conducted synthetic experiments and simulation experiments with real-data on a retail service application.
The results show that the proposed method outputs solutions with higher total revenues than do baselines such as the (modified) retraining method and Bayesian optimization. 

\paragraph{Notation}
Bold lowercase symbols (e.g., $\bx, \by$) denote vectors, and $\|\bx\|$ denotes the Euclidean norm of a vector $\bx$.
The inner product of the vectors $\bx, \by$ is denoted by $\bx^\top \by$.
Let $\mR_+$ be the set of positive real numbers.
The gradient for a real-valued function $f$ \emph{w.r.t.} $\bx$ is denoted by $\nabla_{\bx} f$ and the Jacobian matrix for a vector valued function $\bp$ \emph{w.r.t.} $\bx$ is denoted by $\frac{d\bp}{d\bx}$.
A binomial coefficient of a pair of integers $m$ and $n$ is written as ${}_{m}C_{n}$.
Let $[N]$ be the set of $\{1,2,\dots,N\}$.

\section{Related Works\label{sec:related_work}}

\subsection{Price Optimization Problems with Stochastic Demand}\label{subsec:pricing_problems}
The previous studies on pricing with stochastic demand considered three types of random variable: 
(a) decision-independent random variables included in buyers' purchase behavior;
(b) decision-independent random variables directly included in demand;
(c) decision-dependent random variables included in demand.
Regarding (a), \cite{chawla2010multi} and \cite{correa2017posted} address pricing problems with stochastic behaviors of multiple agents; 
each agent $i$ has a (decision-independent) random variable $v_i$ as its value for a product and purchases a product when the price is below that value.
Regarding (b), \cite{yong2009coordinating,heydari2015two}, and \cite{dong2017electricity} deal with demand with a decision-independent uncertainty, such as $d(x)+\xi$, where $x$ is price and $\xi$ is a random variable independent of price.
Regarding (c), \cite{Bertsimas2005dynamic,wang2019integrating,schulte2020price}, and \cite{hikima2021integrated,hikima2022online,Hikima_Akagi_Kim_Asami_2023}
tackle pricing problems with decision-dependent stochastic demand, such as $d(x)+\xi(x)$, where $\xi(x)$ is a (decision-dependent) random variable.
Our study is categorized into (c).

In this paper, we propose a new general pricing problem with decision-dependent random variables.
Our problem has advantages over the previous ones for tackling (a), (b), and (c).
Regarding (a), while the previous studies need to define agents' actions (e.g., buy or leave) by discontinuous functions, in our formulation, we can define those without a discontinuous function, leading to gradient-based methods.
Regarding (b), we generalize the noise $\xi$ of demand to make it depend on the decision variable, which allows us to deal with situations where the demand uncertainty varies with price.
Regarding (c), previous studies have limited applications since they make application-specific assumptions to obtain an optimal solution or an approximation solution: \cite{wang2019integrating} and \cite{schulte2020price} consider specific situations to optimize prices over multiple periods to sell items and describe efficient methods to find an optimal solution;
\cite{hikima2021integrated,hikima2022online,Hikima_Akagi_Kim_Asami_2023} tackle resource allocation problems while controlling agents' acceptance probabilities for prices and present approximation algorithms with constant approximation ratios;
\cite{Bertsimas2005dynamic} consider a simple demand function where the price of each item does not affect demand for other items and present heuristics to obtain an approximation solution.
In contrast, we deal with a more general framework that has various applications  (see Section \ref{sec:applications}). 
Consequently, our formulation is a non-convex optimization problem and we develop a stochastic method that is theoretically guaranteed to converge to a stationary point.

\subsection{Optimization Methods for Stochastic Problems with Decision-dependent Uncertainty} \label{subsec:related_decision_dependent_noise}
Our price optimization problem, (P) in Section \ref{subsec:problem_p_deginition}, is categorized as a stochastic problem with decision-dependent uncertainty \citep{hellemo2018decision,Varaiya1989}.
This is because the demand of items and services follows a probability distribution depending on price (decision variables). 
Here, we explain three different techniques for solving the problem.\footnote{
Another formulation dealing with decision-dependent random variables is \emph{decision-dependent distributionally robust optimization} \citep{luo2020distributionally,basciftci2021distributionally}.
Although such methods are effective at finding an optimal solution in the worst case when the probability distribution is ambiguous,
they are not appropriate for the purpose of this study.}

\paragraph{Retraining methods \citep{perdomo2020performative,mendler2020stochastic}.}
Retraining methods fix the distribution at each iteration and update the current iterate.
Specifically, \citep{perdomo2020performative} proposed \textit{repeated gradient descent}: $\bx_{k+1}:= \mathrm{proj}_{\mathcal C}(\bx_k - \eta_k \E_{\bmxi \sim D(\bx_k)} [\nabla_{\bx} f(\bx_k,\bmxi)]),$
where $\mathcal C$ is the feasible region and $\mathrm{proj}_{\mathcal C}$ is the Euclidean projection operator onto $\mathcal C$.
It converges to a \emph{performatively stable point} $\bx_{\mathrm{PS}} = \mathrm{arg} \min_{\bx} 
\E_{\bmxi \sim D(\bx_{\mathrm{PS}})}[f(\bx,\bmxi)]$.
However, these methods assume the strong convexity of $f(\bx,\bmxi)$ w.r.t. $\bx$ and are not applicable to our problem.
We provide an intuitive example where RGD fails to work in price optimization, where the objective function is generally not strongly convex.

\begin{example} \label{exa:exisiting_for_pricing}
    Suppose that a seller determines the price $x \in [0, M]$ of a product. 
    The buyer purchases the product ($\xi=1$) with probability $p(x)$ or does not purchase it ($\xi=0$) with probability $1-p(x)$, where $p:[0,M] \to [0,1]$ is a decreasing function.
    The seller wants to solve $\min_{x \in [0,M]} \E_{\xi\sim D(x)}[-x\xi]$ to maximize the expected revenue, where $D(x)$ is the distribution for $\xi$.
    Then, the optimal solution is $x^* \in \arg \min_{x \in [0, M]} -xp(x)$.
    However, RGD continues to raise the price until the purchase probability reaches zero or the price reaches $M$ since $\E_{\xi\sim D(x^k)}[-\nabla_x (x\xi)]=\E_{\xi\sim D(x^k)}[-\xi]=-p(x^k)$ and $p(x^k)\ge 0$ for all $x^k \in [0,M]$.
    This price is generally not equal to $x^*$.
\end{example}

\paragraph{Meta-model methods \citep{brochu2010tutorial,frazier2018tutorial,miller2021outside}.}
This type of method creates a meta-model of the objective function or the distribution map $D(\cdot)$ from multiple sample points.
Bayesian optimization \citep{brochu2010tutorial,frazier2018tutorial} is the process of learning the objective function through Gaussian process regression while finding the global optimal solution.
The two-stage approach \citep{miller2021outside}
estimates a coarse model of the distribution map $D(\cdot)$ and then optimizes a proxy to the objective function by treating the estimated distribution as if it were the true distribution map.
While these methods are powerful for certain problems, they are not suitable for ours:
Bayesian optimization cannot find good solutions when the dimension of the decision variables is too large to be adequately explored;
the two-stage approach assumes that the distribution map is included in location-scale families \citep[Eq. (2)]{miller2021outside}, which cannot be assumed in our problem.

\paragraph{Gradient-free methods \citep{spall2005introduction,Flaxman2004online}.}
Gradient-free methods estimate the gradient by querying objective values at randomly perturbed points around the current iterate.
While this type of method is generic, it often requires many evaluations of objective values to estimate the gradient accurately.

We developed a new projected stochastic gradient descent method by deriving an unbiased stochastic gradient. 
Our method has advantages over the existing ones: unlike retraining methods, it can find stationary points for general pricing problems with no strongly convex objective functions;
unlike meta-model methods, it can find stationary points in high-dimensional optimization problems and does not place a strong assumption on the distribution map; while gradient-free methods naively approximate the gradient, our method approximates it by using gradient information on the objective function and the probability density function, which enables us to estimate gradients more accurately in a shorter computation time.

\section{Optimization Problem \label{sec:problem_definition}}

\subsection{Problem Definition \label{subsec:problem_p_deginition}}
We will consider the following hypothetical situation.
There is a decision maker determining a price vector $\bx \in \mathcal C \subseteq \R^n$ for items $i=1,2,\dots,n$, where the index $i$ denotes the type of items and/or the time period.
Then, the demand vector $\bmxi \in \Xi \subseteq \mR^n$ of $n$ items is sampled from a probability distribution $D(\bx)$.
The decision maker obtains a profit of $s(\bx,\bmxi)-c(\bmxi)$, where $s: \mathcal C \times \Xi \to \R$ and $c:\Xi \to \R$ are the sales and cost functions, respectively.

The revenue maximization problem is as follows:
\begin{align*}
\mathrm{(P)} \quad  \min_{\bx \in \mathcal{C}}&\quad 
\E_{\bmxi \sim D(\bx)}[f(\bx,\bmxi)],
\end{align*}
where $f(\bx,\bmxi):=-s(\bx,\bmxi)+c(\bmxi)$ is real-valued and possibly non-convex.
$D(\bx)$ is a decision-dependent distribution for the measurable set $\Xi$.
Here, we let $\mathrm{Pr}(\bmxi \mid \bx)$ be the probability density function of $D(\bx)$ and assume that the decision maker can obtain the value of $\mathrm{Pr}(\bmxi \mid \bx)$ and $\nabla_{\bx}\mathrm{Pr}(\bmxi \mid \bx)$ for given $\bx$ and $\bmxi$.
This assumption naturally holds in many applications of price optimization.\footnote{For example, 
in \citep{Bertsimas2005dynamic}, the demand for item $i$ at price $x_i$ is defined by $d_i(x_i):=\alpha_i(x_i)+\beta_i(x_i) \xi_i$, where $\xi_i$ is a random variable and its probability density function is given.
In \citep{schulte2020price}, the buyer's arrival rate at price $x$ is assumed to follow a Poisson process with the intensity function $\lambda(x)$, which identifies the probability density function of demand.
\cite{hikima2022online} define $p_{vt}(x)$ as the probability that buyer $v$ arrives at time interval $t$ for price $x$, and then give a probability density function for demand.}

\subsection{Assumptions} \label{sec:assumptions}
Our development of an unbiased stochastic gradient for (P) whose variance is bounded by a constant requires a number of assumptions.
In particular, we will make the following assumptions.
\begin{assumption} \label{assumption: pvt_v2}
For all $\bx \in \mR^n$ and $\bmxi \in \Xi$, the following hold,
\begin{enumerate}
    \item[(i)] $f(\bx,\bmxi)$ is differentiable and Lipschitz continuous with modulus $L_f$ w.r.t. $\bx$ and continuous w.r.t. $\bmxi$,
    \item[(ii)] $\mpr(\bmxi \mid \bx)$ is differentiable w.r.t. $\bx$ and $\mpr(\bmxi \mid \bx)> 0$, and
    \item[(iii)] $\left\|\frac{\nabla_{\bx} \mpr(\bmxi \mid \bx)}{\mpr(\bmxi \mid \bx)} \right\| \le M$ for a constant $M$.
\end{enumerate}
\end{assumption}
\begin{assumption} \label{asp:compact_set}
The set $\mathcal C$ is compact and convex. The set $\Xi$ is compact.
\end{assumption}
Moreover, we need the following assumption when $\bmxi$ is a continuous random vector:
\begin{assumption} \label{asp:continuous_dist}
The set $\Xi$ is a Borel set on $\R^n$. 
Moreover, $\mpr(\bmxi\mid\bx)$ is continuous w.r.t. $\bmxi$ for all $\bx \in \R^n$.
\end{assumption}

Assumptions \ref{assumption: pvt_v2}--\ref{asp:continuous_dist} do not depend on a specific application; there are various applications that satisfy them (see Section \ref{sec:applications}).
Condition (i) of Assumption~\ref{assumption: pvt_v2} usually holds in pricing applications;
the sales function $s(\bx,\bmxi)$ is usually expressed as $\bx^\top \bmxi$ (the product of price and demand), so it can be differentiable w.r.t. $\bx$ and Lipschitz continuous when $\bmxi$ is bounded;
the cost function $c(\bmxi)$ is usually continuous w.r.t. $\bmxi$ since the production cost is usually continuous with respect to demand.
Condition (ii) of Assumption~\ref{assumption: pvt_v2} is satisfied by many distributions with a (statistical) parameter $\blambda(\bx)$, where $\blambda$ is a differentiable vector-valued function.\footnote{For example, the probability density functions of normal and multinomial distributions satisfy condition (ii) of Assumption~\ref{assumption: pvt_v2}.
Since the probability density functions of these distributions are differentiable with respect to their parameters $\blambda$ (e.g., mean, variance), they are also differentiable w.r.t. $\bx$ from the differentiability of $\blambda(\bx)$.}
Condition (iii) of Assumption~\ref{assumption: pvt_v2} means that when the probability of a given demand is small, the effect of price on that probability is also small.
In our application examples presented in Section~\ref{sec:applications}, the multinomial and truncated normal distributions parameterized by price satisfy these conditions.
Assumption \ref{asp:compact_set} is natural for practical pricing applications since price and demand ranges are usually bounded.
Assumption~\ref{asp:continuous_dist} is satisfied if $\bmxi$ follows one of the major continuous probability distributions such as the normal and logistic distributions.
In the next section, we show that our application examples satisfy Assumptions \ref{assumption: pvt_v2}--\ref{asp:continuous_dist}.

\paragraph{Remark.}
Assumption \ref{asp:compact_set} does not hold in the case of unconstrained price optimization, but we can assume $\bx \in [-G,G]^n$ for a sufficiently large $G \in \R_+$ in practice.

\subsection{Application Examples} \label{sec:applications}

\subsubsection{Multiproduct Pricing} \label{subsec:multiple_products}
We consider a variant of \citep{gallego2014multiproduct,zhang2018multiproduct} in which a decision maker exists that determines the prices of multiple products and there are $m$ buyers and $n$ products.
Let $\bx :=(x_1,x_2,\dots,x_n) \in [x_{\min}, x_{\max}]^n$ be the price vector for the products.
We assume buyers choose one product stochastically;
Each buyer chooses product $i \in I:= \{1,\dots,n\}$ with probability
%\begin{align}
$p_i(\bx)= \frac{e^{\gamma_i(\alpha_i-x_i)}}{a_0+\tsum_{j=1}^n e^{\gamma_j(\alpha_j-x_j)}}$ %\label{eq:multinomial_logit}
%\end{align}
or does not choose any product with probability
%\begin{align}
$p_0(\bx)= \frac{a_0}{a_0+\tsum_{j=1}^n e^{\gamma_j(\alpha_j-x_j)}}$.\footnote{Besides the multinomial logit model, various other models can be considered, such as the nested logit model \citep{gallego2014multiproduct} and the generalized nested logit model \citep{zhang2018multiproduct}.}
Here, $\alpha_i$ and $\gamma_i$ are positive constants that can be estimated from historical transaction data \citep{croissant2012estimation}.
Let $\bmxi \in \{ 0,1,\dots,m\}^{n+1}$ be a random vector, where $\xi_0$ represents the number of buyers not purchasing any product and $\xi_i$ for $i=1,\dots,n$ represents the number of sales of each product.
Let $s(\bx,\bmxi)$ and $c(\bmxi)$ be real-valued functions representing the sales and costs of products, respectively.
The following functions are possible for $s$ and $c$:
\small
\begin{align}
&s(\bx,\bmxi):=\sum_{i=1}^n x_i \xi_i, \quad c(\bmxi):= \sum_{i=1}^n c_i(\xi_i), %\label{eq:sell_cost}
&\mathrm{where} \ 
& c_i(\xi_i):=
\begin{cases}
\eta^1_i \xi_i, &\xi_i \le l_i, \\
\eta^2_i (\xi_i -l_i) + \eta^1_i l_i, &l_i < \xi_i \le u_i, \\
\eta^3_i (\xi_i -u_i) + \eta^2_i (u_i -l_i) + \eta^1_i l_i, &\xi_i > u_i.
\end{cases} \nonumber %\label{eq:cost_example_2}
\end{align}
\normalsize
Here, $\eta^1_i$, $\eta^2_i$, $\eta^3_i$, $l_i$, and $u_i$ are constants for each $i$.
The function $c_i$ represents the case where the cost rate varies with the number of sold products (which is also called economies of scale or diseconomies of scale).

The revenue-maximizing problem is as follows:
\begin{align*}
    \min_{\bx\in [x_{\min},x_{\max}]^n} \E_{\bmxi \sim D(\bx)} \left[ -s(\bx,\bmxi) + c(\bmxi) \right],
\end{align*}
where the probability mass function of $D(\bx)$ is $\Pr(\bmxi \mid \bx) := \prod_{i=0}^n {}_{m}C_{\xi_i}
p_i(\bx)^{\xi_i}$.
It can be written in the form of (P).\footnote{If the functions $s$ and $c$ are linear \emph{w.r.t.} $\xi$, then $\E_{\bmxi \sim D(\bx)} \left[ -s(\bx,\bmxi) + c(\bmxi) \right]= -s(\bx,\E_{\bmxi \sim D(\bx)} [\bmxi]) + c(\E_{\bmxi \sim D(\bx)}[\bmxi])$, and the problem is a deterministic optimization as is tackled in \citep{gallego2014multiproduct,zhang2018multiproduct}. Therefore, the problem can be regarded as a generalization of the problems in the previous studies in terms of sales and cost functions.}

The following proposition shows that this application satisfies our assumptions.
\begin{proposition}
Let $\gamma^{\max}:=\max_{i \in I} |\gamma_i|$.
The problem of \textbf{multiproduct pricing} satisfies Assumptions \ref{assumption: pvt_v2} and \ref{asp:compact_set},
where $f(\bx,\bmxi):=-s(\bx,\bmxi) + c(\bmxi)$, $\mathcal C := [x_{\min},x_{\max}]^n$, $L_f:=m$, and $M:= n m \gamma^{\max}$.
\end{proposition}
The proofs of this proposition and the others can be found in Appendix \ref{sec:proof}.

\paragraph{Remark.}
Assumption \ref{assumption: pvt_v2} is not satisfied if $p_i(\bx)$ is non-differentiable.
However, we can satisfy it by using a smoothing technique \citep{chen1996class,chen2012smoothing} to smooth $p_i(\bx)$.

\subsubsection{Congestion Pricing for HOT Lanes} \label{subsec:congestion_pricing}
We consider a stochastic variant of \citep{lou2011optimal} in the following traffic situation:\footnote{Our method can be extended to more general situations, such as ones with many lanes.} there are two lanes, a high-occupancy/toll (HOT) lane and a regular lane; drivers can only switch from the regular lane to the HOT lane.
There is a decision maker determining a price $x_i \in [x_{\min}, x_{\max}]$ of the HOT lane for each time interval $i \in I :=\{1,2,\dots,n\}$.
The purpose of the decision maker is (i) to maximize the total flow rate at the bottlenecks of the HOT and regular lanes and (ii) to prevent the density of vehicles at the switching point from exceeding a certain level (to avoid traffic accidents).
Let $d_i$ be the number of homogeneous drivers in the regular lane in a time interval $i \in I$.
Here, each driver in $i$ changes lane with a probability $p_i(x_i) := \frac{1}{1+e^{\alpha_i h_i+\beta_i x_i+\gamma_i}}$, 
where $h_i$ is a constant indicating the average time savings if a driver chooses the HOT lane at time $i\in I$.
The parameters $\alpha_i$, $\beta_i$, and $\gamma_i$ are constants, which can be estimated in real-time \citep[Section 2.2]{lou2011optimal}.

The optimization problem is as follows (the details can be found in \citep[Section 3]{lou2011optimal}):
\begin{align*}
    \max_{\bx\in [x_{\min},x_{\max}]^{I}} \E_{\bmxi \sim D(\bx)} \left[ \tsum_{i \in I}\left(q_H\left(\xi_i\right)+q_R\left(d_i-\xi_i\right)\right) + \theta \min \left(\tilde{k}-\frac{1}{|I|} \tsum_{i \in I} k\left(\xi_i\right),0\right) \right],
\end{align*}
where $\xi_i \in \{0,1,\dots, d_i\}$ is a random variable indicating the number of drivers switching their lanes in $i$.
Regarding the first term, the values of $q_H\left(\xi_{i}\right)$ and $q_R\left(d_i-\xi_{i}\right)$ are continuous functions representing flow rates at the bottlenecks on the HOT lane and the regular lane, respectively. 
This term aims to maximize the flow rate of each lane.
Regarding the second term, $\tilde{k}\in \mR_+$ is the critical density (the density likely to cause traffic accidents) of vehicles at the switching point, and $k(\xi_i)$ is a continuous function representing the density at the switching point for the demand $\xi_i$ in $i$.
Therefore, $\theta \min \left(\tilde{k}-\frac{1}{|I|} \tsum_{i \in I} k\left(\xi_i\right),0\right)$ is a penalty term for densities above the critical density, where $\theta \ge 0$ is the penalty parameter.
This optimization problem can be written in the form of (P), where $s(\bx,\bmxi):=0$, $c(\bmxi):=\tsum_{i \in I}\left(q_H\left(\xi_i\right)+q_R\left(d_i-\xi_i\right)\right) + \theta \min \left(\tilde{k}-\frac{1}{|I|} \tsum_{i \in I} k\left(\xi_i\right),0\right)$, and
$\Pr(\bmxi \mid \bx) :=\prod_{i \in I} 
{}_{d_i}C_{\xi_i}
p_i(x_i)^{\xi_i}(1-p_i(x_i))^{d_i-\xi_i}$.

The following proposition shows that this application satisfies our assumptions.
\begin{proposition}
The problem of \textbf{congestion pricing for HOT lanes} satisfies Assumptions \ref{assumption: pvt_v2} and \ref{asp:compact_set},
where $f(\bx, \bmxi):= \tsum_{i \in I}\left(q_H\left(\xi_i\right)+q_R\left(d_i-\xi_i\right)\right)+ \theta \min (\tilde{k}-\frac{1}{|I|} \tsum_{i \in I} k\left(\xi_i\right),0)$,
    $L_f:=0$, $\mathcal C := [x_{\min},x_{\max}]^{I}$, and $M:= |I| \max_{i \in I}\left( |\beta_i| d_i \right)$.
\end{proposition}

\subsubsection{Pricing with Demand Prediction from Limited Data Points} \label{subsec:pricing_limited_data}
Here, we will consider optimizing prices of $n$ types of item.
Regarding the prices $\bx \in [x_{\min},x_{\max}]^n$,
the demand $\xi_i\in \{\xi \mid 0 \le \xi\le \xi_i^{\max}\}$ of item $i$ is predicted using data points $\hat{D}:=\{(\hat{\bx}^d, \hat{\bmxi}^d)\}_{d=1}^N$ through the truncated Gaussian process \citep[Section 8.1]{swiler2020survey}:
\begin{align*}
    &\xi_i \sim \frac{\bm{1}_{\{0\le \xi_i\le\xi_i^{\max}\}}(\xi_i)}{C^i(\bx)} 
    N[\bv^i(\bx)^\top \ba^i, %\bm{K}^{-1}\by,
    (\sigma^i)^2-\bv^i(\bx)^\top %\bm{K}^{-1} \bv^i(\bx)],
    A^i \bv^i(\bx)], \\
    &\textrm{where} \  C^i(\bx)=\int_{0}^{\xi_i^{\max}} \frac{1}{\sqrt{2\pi ((\sigma^i)^2-\bv^i(\bx)^\top A^i \bv^i(\bx))}} \exp\left(-\frac{(\phi-\bv^i(\bx)^\top \ba^i)^2}{2((\sigma^i)^2-\bv^i(\bx)^\top A^i \bv^i(\bx))}\right) d\phi.
\end{align*}
Here, $\bv^i: \R^n \to \R^N$, $\ba^i \in \R^N$, $\sigma^i\in \R$, and $A^i\in \R^{N\times N}$ are respectively the function, vector, scalar, and matrix learned from the data points $\hat{D}$.
The $j$-th element $v_j^i(\bx)$ of $\bv^i(\bx)$ is defined by $
v^i_j(\bx):=\theta_1^i \exp \left(-\frac{\|\bx-\hat{\bx}^j\|^2}{\theta^i_2}\right)$, 
%\label{eq:gaussian_kernel}
%\end{align}
where $\theta_1^i\in\R_+$ and $\theta_2^i\in\R_+$ are learned constants.
The normalization function $C^i(\bx)$ is the probability that a sample lies in $\{0\le \xi_i \le \xi_i^{\max}\}$.
Here, $(\sigma^i)^2-\bv^i(\bx)^\top A^i \bv^i(\bx) \ge \Delta$ for some $\Delta \in \R_+$.\footnote{Given that the observations are subject to noise, it is natural to predict that the variance $(\sigma^i)^2-\bv^i(\bx)^\top A^i \bv^i(\bx)$ is more than or equal to a certain constant ($\Delta$).}

The revenue-maximizing problem is as follows:
\begin{align*}
    \min_{\bx \in [x_{\min},x_{\max}]^n} \ \E_{\bmxi \sim D(\bx)} \left[-s(\bx,\bmxi) + c\left(\bmxi   \right) \right],
\end{align*}
where $s(\bx,\bmxi):= \bmxi^\top \bx$, $c(\bmxi):= \sum_{i=1}^n c_i(\xi_i)$, and $c_i:\R \to \R$ is a continuous function for $i=1,2,\dots,n$.
$c_i(\xi_i)$ represents the cost for item $i$.
This problem can be written in the form of (P), where
%\begin{align*}
 $\Pr(\bmxi \mid \bx):= 
\prod_{i=1}^n \frac{1}{C^i(\bx)\sqrt{2\pi((\sigma^i)^2-\bv^i(\bx)^\top A^i \bv^i(\bx))}} \exp\left(-\frac{(\xi_i-\bv^i(\bx)^\top \ba^i)^2}{2((\sigma^i)^2-\bv^i(\bx)^\top A^i \bv^i(\bx))}\right)   $
%\end{align*}
for $\bmxi \in \Xi := \{\bmxi \mid \forall i \in [n], 0 \le \xi_i \le \xi_i^{\max} \}$.

The following proposition shows that this application satisfies our assumptions.
\begin{proposition}
Let
$\theta_1^{\max}:=\max_{i} \theta_1^i$, $\theta_2^{\min}:=\min_{i} \theta_2^i$, $a^{\max}:=\max_{i,k} |a^i_k|$, $A^{\max}:=\max_{i,k,l} |A^i_{k,l}|$, and $\xi^{\max}:=\max_{i} \xi_i^{\max}$.
The problem of \textbf{pricing with demand prediction from limited data points}  satisfies Assumptions \ref{assumption: pvt_v2}--\ref{asp:continuous_dist},
where $f(\bx,\bmxi):=-s(\bx,\bmxi) + c(\bmxi)$, $\mathcal C := [x_{\min},x_{\max}]^n$, $L_f:=n\xi^{\max}$, and $M:= \frac{4n^2N\theta_1^{\max}(x_{\max}-x_{\min})}{\Delta \theta_2^{\min}} \left(NA^{\max}\theta_1^{\max} + (\xi^{\max} + N\theta_1^{\max} a^{\max}) \left(a^{\max}+NA^{\max}\theta_1^{\max}\frac{\xi^{\max} + N\theta_1^{\max} a^{\max}}{\Delta}\right) \right)$.
\end{proposition}

\subsection{Advantages of Our Formulation \label{sec:importance}}
\subsubsection{Benefits of Using Decision-dependent Random Variables} \label{subsec:related_independent_formlulation}
The \textbf{multiproduct pricing} problem in Section \ref{subsec:multiple_products} can also be expressed in terms of decision-independent random variables as follows.
Each buyer $j=1,2,\dots,m$ has a value $\gamma_i(\alpha_i-x_i)+\mu_{ij}$ for each product $i$, where $\alpha_{i}$ and $\gamma_i$ are constants, $x_i$ is the price, and $\mu_{ij}$ is a random variable following a Gumbel distribution with mode 0 and variance $\frac{\pi^2}{6}$.
Each buyer purchases a product $i$ with the highest $\gamma_i(\alpha_i-x_i)+\mu_{ij}$.
Accordingly, the demand $\xi_i$ for product $i$ can be defined by
$\xi_i(\bx,\bmu):=\sum_{j =1}^m \xi_{ij}(\bx,\bmu_j)$,
where $\bmu_j =(\mu_{ij})_{i=1}^n$ and $\xi_{ij}(\bx,\bmu_j):=1$ if $i= \argmax_{r} \{ \gamma_r(\alpha_r-x_r)+\mu_{rj}\}$ with $\xi_{ij}(\bx,\bmu_j):=0$ otherwise.
The optimization problem can be written as follows by letting $\bmxi(\bx,\bmu) := (\xi_i(\bx,\bmu))_{i=1}^n$:
\begin{align*}
    \min_{\bx\in [x_{\min},x_{\max}]^n} \E_{\bmu} \left[ -s(\bx,\bmxi(\bx,\bmu)) + c(\bmxi(\bx,\bmu)) \right].
\end{align*}

Although the \textbf{multiproduct pricing} problem can be formulated in the above manner  with decision-independent random variables, the discontinuous function $\xi_{ij}(\bx,\bmu_j)$ makes it difficult to optimize.\footnote{Optimization problems involving such discontinuous functions have been addressed by \citep{correa2017posted}. 
They propose approximation methods to deal with this difficulty.}
In contrast, our problem does not involve a discontinuous function,  which allows us to use gradient-based methods.

Moreover, \cite{yong2009coordinating} and \cite{dong2017electricity} tackle similar problems to ours by defining
the demand for a product/service as $d(x)+\xi$, where $x$ is price and $\xi$ is a decision-independent random variable. 
However, assuming $\xi$ is decision-independent makes it impossible to handle situations where demand uncertainty varies with price.
In contrast, our problem setting can deal with such a situation by using decision-dependent random variables.

\subsubsection{Differences from Existing Pricing Problems with Decision-dependent Uncertainty} \label{subsec:def_existing_dd}
The existing formulations with decision-dependent uncertainty make assumptions specific to their applications.
For example, \cite{schulte2020price} assume that demand follows a Poisson distribution with an intensity $\lambda(x)$, and they cannot use a multinomial or truncated Gaussian distribution as the demand distribution.
Moreover, since a fixed cost is charged on their products, they can not handle a nonlinear cost.
\cite{hikima2021integrated} assume that
the probability density function for demand is ${\rm Pr}(\bmxi \mid \bx)=\displaystyle{ \prod\nolimits_{u\in U}} \left\{ p_u(x_u)^{\xi_u}(1-p_u(x_u))^{(1-\xi_u)}  \right\}$, where $p_u(x)$ is the probability that service user $u \in U$ accepts the price $x_u$; they cannot use a truncated Gaussian distribution.
In addition, they assume a specific objective function, which is defined by a bipartite matching problem with uncertainty.

In contrast to the existing formulations, ours has more varied applications because it has more general assumptions.
The trade off for this generality is that our problem is non-convex and the dependence of the probability distribution on the decision variables defeats conventional stochastic optimization theory.
Below, we focus on finding a stationary point and develop a projected stochastic gradient descent method by deriving unbiased stochastic gradient estimators.

\section{Proposed Method \label{sec:proposed_method}}

\subsection{Preliminaries}

\begin{definition}[Projection oracle] 
Given a point $\bx$, we define the following as a projection oracle:
\begin{align}
\mathrm{proj}_{\mathcal C}(\bx):= \mathrm{arg} \min_{\by \in \R^n} \{ \|\bx -\by \|_2 \mid \by \in \mathcal C \}.\nonumber
\end{align}
\end{definition}

\begin{definition}[Unbiased stochastic gradient] \label{USG}
Given a point $\bx$, we call $g(\bx,\bmxi)$ an ``unbiased stochastic gradient" if $$\E_{\bmxi \sim D(\bx)}[g(\bx,\bmxi)]=\nabla_{\bx} \E_{\bmxi \sim D(\bx)}[f(\bx,\bmxi)].$$
\end{definition}

\begin{definition}[Gradient mapping] \label{def:gradinet_mapping}
Given a point $\bx$ and $\eta \in \mR_+$, the gradient mapping of (P) is defined by
$$\mathcal G(\bx,\eta) := \frac{1}{\eta}(\bx - \mathrm{proj}_{\mathcal C}(\bx-\eta \nabla_{\bx} \E_{\bmxi \sim D(\bx)}[f(\bx,\bmxi)])).$$
\end{definition}
\begin{definition}[$\epsilon$-stationary point] \label{def:epsilon_statinary}
We call $\hat{\bx}$ an $\epsilon$-stationary point for (P) if $ \E_{\hat{\bx}} [\|\mathcal G (\hat{\bx},\eta)\|^2] \le \epsilon^2$ for some $\eta \in \mR_+$, where $\hat{\bx}$ denotes the point returned by a stochastic algorithm.
\end{definition}

The following preliminary lemmas are needed for ensuring our method's convergence when the random variables are continuous.
\begin{lemma} \label{lem:int_lim_change}
Suppose that Assumptions \ref{assumption: pvt_v2}--\ref{asp:continuous_dist} hold.
Let $h(\bx,\bmxi):=f(\bx,\bmxi) \mpr(\bmxi\mid\bx)$.
Then,  
\begin{align*}
\nabla_{\bx} \int_{\bmxi \in \Xi} h(\bx,\bmxi) d\bmxi = \int_{\bmxi \in \Xi} \nabla_{\bx} h(\bx,\bmxi)d{\bmxi}
\end{align*}
for all $\bx \in \mathcal C$.
\end{lemma}

\begin{lemma} \label{lem:int_lim_change_2}
Suppose that conditions (ii) and (iii) of Assumption \ref{assumption: pvt_v2} and Assumptions \ref{asp:compact_set} and \ref{asp:continuous_dist} hold.
Then, 
\begin{align*}
  \nabla_{\bx} \int_{\bmxi \in \Xi} q(\bmxi)\mpr(\bmxi\mid\bx) d\bmxi = \int_{\bmxi \in \Xi} \nabla_{\bx} q(\bmxi)\mpr(\bmxi\mid\bx)d{\bmxi}
\end{align*}
for all $\bx \in \mathcal C$
  and any real-valued continuous function $q:\Xi \to \R$.
\end{lemma}

Throughout the paper, we let 
\begin{align}
f_{\max}:=\max_{\bx \in \mathcal C, \bmxi \in \Xi} |f(\bx,\bmxi)|, \label{def:f_max}
\end{align}
which exists since $\mathcal C$ and $\Xi$ are compact from Assumption \ref{asp:compact_set} and $f(\bx,\bmxi)$ is real-valued and continuous from Assumption \ref{assumption: pvt_v2}.

\subsection{Unbiased Stochastic Gradient for (P)}
First, we propose an unbiased stochastic gradient for (P). 

\begin{lemma} \label{sfo_lemma}
Suppose that conditions (i) and (ii) of Assumption~\ref{assumption: pvt_v2} hold.
Moreover, suppose that condition (iii) of Assumption \ref{assumption: pvt_v2}, Assumption \ref{asp:compact_set}, and Assumption \ref{asp:continuous_dist} hold if $\bmxi$ is a continuous random vector.
Let $\delta \in \R$ and
$$\bg(\bx,\bmxi,\delta):=\nabla_{\bx} f(\bx,\bmxi)+\left( f(\bx,\bmxi) - \delta \right) \frac{\nabla_{\bx} \mpr(\bmxi\mid\bx)}{\mpr(\bmxi\mid\bx)}.$$
Then, $\bg(\bx,\bmxi,\delta)$ is an unbiased stochastic gradient for (P) for any $\delta \in \R$.
\end{lemma}
Inspired by a technique called \textit{baseline} in reinforcement learning \citep{williams1992simple,sutton2018reinforcement}, we decided to include a variance reduction parameter $\delta$ in the unbiased stochastic gradient. 
If $\delta$ is close to $\E_{\bmxi \sim D(\bx)}[f(\bx,\bmxi)]$, the second term of $g(\bx,\bmxi,\delta)$ is small, and the variance of $g(\bx,\bmxi,\delta)$ is reduced. 
We show how to determine $\delta$ in Section \ref{sec:delta_detemination}.

The gradient in Lemma~\ref{sfo_lemma} has the following useful feature.
\begin{lemma} \label{lemma:SFO_variance}
Suppose that Assumptions \ref{assumption: pvt_v2} and \ref{asp:compact_set} hold.
Moreover, suppose that Assumption \ref{asp:continuous_dist} holds if $\bmxi$ is a continuous random vector.
Let $\delta \in [-f_{\max}, f_{\max}]$.
Then, for all $\bx \in \mathcal C$,
\begin{align*}
\E_{\bmxi' \sim D(\bx)}[\| \bg(\bx,\bmxi',\delta) - \nabla_{\bx} \E_{\bmxi \sim D(\bx)}[f(\bx,\bmxi)]  \|^2] 
\le (L_f+2f_{\max}M)^2,
\end{align*}
where $\bg(\bx,\bmxi',\delta)$ is defined as in Lemma \ref{sfo_lemma}.
\end{lemma}
Lemma \ref{lemma:SFO_variance} shows that the variance of the stochastic gradient of Lemma~\ref{sfo_lemma} can be bounded by a constant.
This is a necessary condition for stochastic gradient methods to have a convergence rate independent of the number of possible values of $\bmxi$ \citep{li2018simple}.

Moreover, the following lemma is necessary for ensuring the convergence of the proposed method.
\begin{lemma} \label{lem:entire_derivative_bound}
Suppose that Assumptions \ref{assumption: pvt_v2} and \ref{asp:compact_set} hold.
Moreover, suppose that Assumption \ref{asp:continuous_dist} holds if $\bmxi$ is a continuous random vector.
Then, for all $\bx \in \mathcal C$,
\begin{align*}
\|\nabla_{\bx} \E_{\bmxi \sim D(\bx)}[f(\bx,\bmxi)]\|\le L_f + f_{\max} M.
\end{align*}
\end{lemma}

\subsection{Calculation of Variance Reduction Parameter $\delta$} \label{sec:delta_detemination}
To reduce the variance of the gradient in Lemma \ref{sfo_lemma}, the parameter $\delta$ should be close to $\E_{\bmxi \sim D(\bx_k)}[f(\bx_k,\bmxi)]$ for the iterate $\bx_k$.
During the iterations of the algorithm, $\delta_k$ is updated to bring it closer to the target value. We consider the following sequential stochastic problem for $\delta_k$:
a decision maker selects $\delta_k\in [-f_{\max}-\kappa, f_{\max}+\kappa]$ at iteration $k$ and incurs an unobserved cost $\psi_k(\delta_k):=\frac{1}{2}(\delta_k-\E_{\bmxi \sim D(\bx_k)}[f(\bx_k,\bmxi)])^2$, where $\bx_k$ and $\kappa \in \R_+$ are given;\footnote{$\kappa$ is an arbitrarily small positive value.
It extends the range of $\delta_k$, which is needed in Proposition \ref{prop:ogd_regret}.}
for the decision maker, an unbiased estimate of $\E_{\bmxi \sim D(\bx_k)}[f(\bx_k,\bmxi)]$, denoted by $v_k$, is obtained by sampling.
Here, we assume $v_k \in [-f_{\max},f_{\max}]$, which is usually holds from the definition \eqref{def:f_max} of $f_{\max}$.

As a way to solve the above problem, we propose Algorithm \ref{alg:OGD}, which is based on the online gradient descent (OGD) algorithm \citep{besbes2015non}.
\begin{algorithm}[t!]
  \caption{OGD algorithm}
  \label{alg:OGD}
  \begin{algorithmic}[1]
    \Require{iteration limit $R$, initial parameter $\delta_1 \in [-f_{\max}, f_{\max}]$, stepsize parameter $\{\zeta_k\}_{k=2}^R$, noisy target value $\{v_k\}_{k\in[R]}$}.
    \Ensure{$\delta_k$ for $k\in[R]$}
    \For{$k =1,\dots, R-1$}
    \State $\delta_{k+1}\gets (1-\zeta_{k+1})\delta_{k}+\zeta_{k+1} v_k$ \label{state:delta_update}
    \EndFor
  \end{algorithmic}
\end{algorithm}
This method updates $\delta_k$ by using the stochastic gradient since $(1-\zeta_{k+1})\delta_{k}+\zeta_{k+1} v_k= \delta_{k} - \zeta_{k+1} (\delta_{k} - v_k) \approx \delta_{k} - \zeta_{k+1} \nabla \psi_k(\delta_{k})$. 
Note that $\nabla \psi_k(\delta_k)=\delta_k-\E_{\bmxi \sim D(\bx_k)}[f(\bx_k,\bmxi)]$ from the definition of $\psi_k$.
Accordingly, the following proposition holds from \citep[Lemma C-5]{besbes2015non}, which guarantees that Algorithm \ref{alg:OGD} outputs a solution close to the optimum in terms of regret.
\begin{proposition} \label{prop:ogd_regret}
Let $\zeta_k:=\frac{1}{k}$ for $k\in [R]$, and let $\delta_k$ be the output of Algorithm \ref{alg:OGD} for $k \in[R]$.
Then, there exists a constant $\bar{C}$ such that
\begin{align*}
    \E\left[\sum_{k=1}^R \psi_k(\delta_k)\right]-\min_{\delta} \sum_{k=1}^R \psi_k(\delta) \le \bar{C} \log R.
\end{align*}
\end{proposition}
From this proposition and the definition of $\psi_k$, we find that output $\delta_k$ of Algorithm \ref{alg:OGD} is a reasonable approximation of $\E_{\bmxi \sim D(\bx_k)}[f(\bx_k,\bmxi)]$.

\subsection{Proposed Algorithm}
We propose Algorithm~\ref{alg:rsag} for solving problem (P).
It incorporates our stochastic gradient and Algorithm~\ref{alg:OGD} into a projected stochastic gradient method \citep[Algorithm 4]{ghadimi2016accelerated}.
Lines \ref{line:x_md_set}--\ref{projection2} update the iterate on the basis of \citep[Algorithm 4]{ghadimi2016accelerated} by using our proposed stochastic gradient.
Line~\ref{line:calc_baseline} updates the variance reduction parameter on the basis of Algorithm \ref{alg:OGD} by letting $v_k$ be $\frac{1}{m_k} \sum_{\ell=1}^{m_{k}}f(\bx_{k}^{md},\bmxi^\ell)$.
Note that line \ref{line:calc_baseline} does not impose any additional computation cost since $\frac{1}{m_k} {\sum_{\ell=1}^{m_{k}}}f(\bx_{k}^{md},\bmxi^\ell)$ is already computed on line~\ref{g_k_generate}.

Then, from Lemmas~\ref{sfo_lemma}--\ref{lem:entire_derivative_bound} and \citep[Corollary 6]{ghadimi2016accelerated}, the following convergence theorem holds.

\begin{algorithm}[t!]
  \caption{Projected stochastic gradient algorithm}
  \label{alg:rsag}
  \begin{algorithmic}[1]
    \Require{initial iterate $\bx_0 \in \mathcal C$, initial variance reduction parameter $\delta_1 \in [-f_{\max},f_{\max}]$, iteration limit $\{N\ge 1\}$, batch size $\{m_k\}_{k\in[N]}$, probability distribution $D_R(N)$ for set $\{1,2,\dots,N\}$, and stepsize parameters $\{\alpha_k \in (0,1]\}_{k\in[N]}$, $\{ \beta_k \in \mR_+\}_{k\in[N]}$, $\{\lambda_k \in \mR_+ \}_{k\in[N]}$}, and $\{\zeta_k\}_{k=2}^N$
    \Ensure{$\bx^{md}_R$}
    \State Set $\bx^{ag}_0=\bx_0$ and sample $R\sim D_R(N)$.
    \For{$k = 1,2,\dots,R$}
    \State $\bx_k^{md} \gets (1-\alpha_k)\bx_{k-1}^{ag}+\alpha_k \bx_{k-1}$  \label{line:x_md_set}
    \State Sample $\bmxi^{\ell} \sim D(\bx_k^{md})$ for $\ell =1,\dots,m_k$.
    \State $\bg_k \gets \frac{1}{m_k} {\displaystyle \sum_{\ell=1}^{m_k}}\left(\nabla_{\bx} f(\bx_k^{md},\bmxi^\ell)  +(f(\bx_k^{md},\bmxi^\ell)-\delta_{k}) \frac{\nabla_{\bx} \mpr(\bmxi^\ell \mid \bx_k^{md})}{\mpr(\bmxi^\ell \mid \bx_k^{md})}\right)$\label{g_k_generate}
     \State $\bx_k\gets \mathrm{proj}_{\mathcal C}(\bx_{k-1}-\lambda_k \bg_k)$ \label{projection1}
        \State $\bx_k^{ag}\gets \mathrm{proj}_{\mathcal C}(\bx_k^{md}-\beta_k \bg_k)$\label{projection2}
    \If{$k \le R-1$} $\delta_{k+1} \gets  (1-\zeta_{k+1}) \delta_{k} +  \frac{\zeta_{k+1}}{m_k} {\displaystyle \sum_{\ell=1}^{m_{k}}}f(\bx_{k}^{md},\bmxi^\ell)$ \label{line:calc_baseline}
    \EndIf
    \EndFor
  \end{algorithmic}
\end{algorithm}

\begin{theorem} \label{thm:convergence}
Suppose that Assumptions \ref{assumption: pvt_v2} and \ref{asp:compact_set} hold.
Moreover, suppose that Assumption \ref{asp:continuous_dist} holds if $\bmxi$ is a continuous random vector.
Let the inputs of Algorithm \ref{alg:rsag} be
$\alpha_k:=\frac{2}{k+1},\ 
 \beta_k:=\frac{1}{2L_{Ef}},\  
 \lambda_k:=\frac{k\beta_k}{2}, m_k:=\left\lceil \frac{(L_f + 2f_{\max} M)^2k}{L_{Ef}\tilde{D}^2}\right\rceil$, and 
$\Pr(R=k):=\frac{\Gamma_k^{-1}\beta_k(1-L_{Ef}\beta_k)}{\sum_{{\tau}=1}^N \Gamma_{\tau}^{-1} \beta_{\tau} (1-L_{Ef}\beta_{\tau})}$ for $k=1,2,\dots,N,$
   where $\tilde{D}$ is some parameter, $L_{Ef}:= L_f + f_{\max} M$, $\Gamma_1:=1$, and $\Gamma_k:=(1-\alpha_k)\Gamma_{k-1}$ for $k=2,\dots,N$.
Let $\zeta_k:=\frac{1}{k}$ for $k=2,\dots,N$.
Then,
\begin{align*}
    \E[\| \mathcal{G}(\bx_R^{md},{\beta_R}) \|^2] 
    \le 96 L_{Ef}\left[ \frac{4L_{Ef} \|\bx_0-\bx^* \|^2}{N(N+1)(N+2)} +\frac{L_{Ef}(\|\bx^*\|^2+H^2)+2\tilde{D}^2}{N}\right],
\end{align*}
where $H:=\max_{\bx\in \mathcal C} \|\bx\|$.
Consequently, to obtain an $\epsilon$-stationary point of Definition \ref{def:epsilon_statinary},
we need at most $O\left(\left[ \frac{L_{Ef}^2\|\bx_0-\bx^*\|^2}{\epsilon^2} \right]^{\frac{1}{3}} + \frac{L_{Ef}^2(\|\bx^*\|^2+H^2)+L_{Ef}\tilde{D}^2}{\epsilon^2} \right)$ iterations.
\end{theorem}
The parameter $\tilde{D}$ in Theorem \ref{thm:convergence} determines the balance between the minibatch size and the iteration complexity: a small $\tilde{D}$ results in smaller iteration complexity but a larger minibatch size; a large $\tilde{D}$ leads to a smaller minibatch size but larger iteration complexity.

\paragraph{Bottleneck of Algorithm \ref{alg:rsag}.}
The bottleneck of Algorithm~\ref{alg:rsag} is line~\ref{g_k_generate} because it requires $\frac{\nabla_{\bx} \mpr(\bmxi^\ell \mid \bx_k^{md})}{\mpr(\bmxi^\ell \mid \bx_k^{md})}$ to be computed $m_k$ times. 
This calculation takes a lot of time since $m_k$ has to be at least proportional to the number $k$ of current iterations to obtain the convergence rate of Theorem~\ref{thm:convergence}.

\subsection{Specialized Projected Stochastic Gradient Method for Price Optimization in Multi-agent Applications}
\label{sec:specializedSGD}
To reduce the computation cost at the bottleneck of Algorithm \ref{alg:rsag}, we propose a specialized projected stochastic gradient method that adds the following assumptions to (P).

\begin{assumption} \label{asp:revenue_average}
$\E_{\bmxi \sim D(\bx)}[s(\bx,\bmxi)] = s(\bx,\E_{\bmxi \sim D(\bx)}[\bmxi])$ and $c(\bmxi)$ is continuous.
\end{assumption}

\begin{assumption} \label{asp:choice_model}
The probability density function $\Pr(\bmxi \mid \bx)$ is defined as $\phi(\bp(\bx),\bmxi)$, where $\phi$ is real-valued and differentiable \emph{w.r.t.} $\bp$, $\frac{\nabla_{\bp} \phi(\bp,\bmxi)}{\phi(\bp,\bmxi)}$ is easily computed, and $\bp$ is vector-valued and differentiable.
\end{assumption}

The above assumptions are often satisfied in price optimization for multi-agent applications. 
In particular, Assumption \ref{asp:revenue_average} tends to hold because the sales function $s$ in price optimization is usually linear with respect to $\bmxi$ and the cost function $c$ is usually continuous with respect to $\bmxi$.
Assumption \ref{asp:choice_model} holds for many parameterized distributions (e.g., binomial, multinomial, and Poisson distributions) since the probability density function and its gradient can be simply written by its parameters.
Many multi-agent applications satisfy Assumption \ref{asp:choice_model} because the distribution of the demand follows a binomial or multinomial distribution with parameters $\bp$, which represents the probabilities of each agent's actions.

The following lemmas show that the applications with multiple agents described in Section \ref{sec:applications} satisfy Assumptions \ref{asp:revenue_average} and~\ref{asp:choice_model}.

\begin{proposition} \label{lem:mpp_ass23}
The problem of \textbf{multiproduct pricing} satisfies Assumption \ref{asp:revenue_average}.
Moreover, let $\phi(\bp(\bx),\bmxi):= \prod_{i=0}^n {}_{m}C_{\xi_i}
p_i(\bx)^{\xi_i}.$
Then, $\Pr(\bmxi \mid \bx)=\phi(\bp(\bx),\bmxi)$ and 
$\left(\frac{\nabla_{\bp}\phi(\bp(\bx),\bmxi)}{\phi(\bp(\bx),\bmxi)}\right)_k=\frac{\xi_k}{p_k(\bx)}.$
\end{proposition}

\begin{proposition} \label{lem:cphl_asp23}
The problem of \textbf{congestion pricing for Hot lanes} satisfies Assumption \ref{asp:revenue_average}.
Moreover, let $\phi(\bp(\bx),\bmxi):=\prod_{i \in I} {}_{d_i}C_{\xi_i}
p_i(x_i)^{\xi_i}(1-p_i(x_i))^{d_i-\xi_i}.$
Then, $\Pr(\bmxi \mid \bx)=\phi(\bp(\bx),\bmxi)$ and 
$\left(\frac{\nabla_{\bp}\phi(\bp(\bx),\bmxi)}{\phi(\bp(\bx),\bmxi)}\right)_k=\frac{\xi_k}{p_k(x_k)}-\frac{d_k-\xi_k}{1-p_k(x_k)}.$
\end{proposition}

Below, we present the lemmas for our specialized method under Assumptions \ref{assumption: pvt_v2}--\ref{asp:choice_model}.
Let $c_{\max}:=\max_{\bmxi \in \Xi} |c(\bmxi)|$, which exists since $\Xi$ is compact from Assumption \ref{asp:compact_set} and $c(\bmxi)$ is continuous from Assumption~\ref{asp:revenue_average}.

\setcounter{lemma}{12}
\begin{lemma} \label{lem:sfo_accel}
Suppose that condition (ii) of Assumption~\ref{assumption: pvt_v2} and Assumptions \ref{asp:revenue_average} and \ref{asp:choice_model} hold.
Moreover, suppose that condition (iii) of Assumption~\ref{assumption: pvt_v2}, Assumption \ref{asp:compact_set}, and Assumption \ref{asp:continuous_dist} hold if $\bmxi$ is a continuous random vector.
Let 
%\begin{align*}
$\bg_2(\bx,\bmxi',\delta):=-\nabla_{\bx} s(\bx,\E_{\bmxi \sim D(\bx)}[\bmxi]) 
+ (c(\bmxi')-\delta) \frac{d \bp(\bx)}{d \bx} \frac{\nabla_{\bp} \phi(\bp(\bx),\bmxi')}{\phi(\bp(\bx),\bmxi')}.$
%\end{align*}
%where $\bmxi^k$ is a random vector following the probability distribution $D(\bx)$.
Then, $\bg_2(\bx,\bmxi',\delta)$ is an unbiased stochastic gradient for (P).
\end{lemma}

\begin{lemma} \label{lemma:SFO_variance_2}
Suppose that 
conditions (ii) and (iii) of Assumption \ref{assumption: pvt_v2}, and Assumptions \ref{asp:compact_set}, \ref{asp:revenue_average}, and \ref{asp:choice_model} hold.
Moreover, suppose that Assumption \ref{asp:continuous_dist} holds if $\bmxi$ is a continuous random vector.
Let $\delta \in [-c_{\max}, c_{\max}]$.
Then, for all $\bx \in \mathcal C$,
%\begin{align*}
$\E_{\bmxi' \sim D(\bx)}[\| \bg_2(\bx,\bmxi',\delta) - \nabla_{\bx} \E_{\bmxi \sim D(\bx)}[f(\bx,\bmxi)]\|^2] \le 4 \left(c_{\max} M\right)^2,$
%\end{align*}
where $\bg_2(\bx,\bmxi',\delta)$ is defined as in Lemma~\ref{lem:sfo_accel}.
\end{lemma}

Now, let us examine Algorithm~\ref{alg:rsag_specialized}. Its computational cost is lower than that of Algorithm \ref{alg:rsag}: Algorithm~\ref{alg:rsag} requires $m_k$ calculations of $\frac{\nabla_{\bx} \mpr(\bmxi^\ell \mid \bx_k^{md})}{\mpr(\bmxi^\ell \mid \bx_k^{md})}$, whereas Algorithm~\ref{alg:rsag_specialized} requires $m_k$ calculations of $\frac{\nabla_{\bp} \phi(\bp(\bx_k^{md}),\bmxi^{\ell})}{\phi(\bp(\bx_k^{md}),\bmxi^{\ell})}$, which can be easily computed from Assumption~\ref{asp:choice_model}.

\begin{algorithm}[t!]
  \caption{Specialized projected stochastic gradient algorithm}
  \label{alg:rsag_specialized}
  \begin{algorithmic}[]
\State In Algorithm \ref{alg:rsag}, let $\delta_1 \in [-c_{\max},c_{\max}]$ in the input, replace line \ref{g_k_generate} by
\begin{flalign*}
\bg_k \gets -\nabla_{\bx} s(\bx_k^{md},\E_{\bmxi \sim D(\bx_k^{md})}[\bmxi]) 
    + \frac{d\bp(\bx_k^{md})}{d\bx} \frac{1}{m_k} \sum_{\ell=1}^{m_k} \left(c(\bmxi^{\ell})-\delta_k\right)\frac{\nabla_{\bp} \phi(\bp(\bx_k^{md}),\bmxi^{\ell})}{\phi(\bp(\bx_k^{md}),\bmxi^{\ell})},&&
\end{flalign*}
\State and replace line~\ref{line:calc_baseline} by
\If{$k \le R-1$} $\delta_{k+1} \gets  (1-\zeta_{k+1}) \delta_{k} +  \frac{\zeta_{k+1}}{m_k} {\displaystyle \sum_{\ell=1}^{m_{k}}}c(\bmxi^\ell)$
\EndIf
%\begin{align*}
%&\delta_{k+1} \gets  (1-\zeta_{k+1}) \delta_{k} +  \frac{\zeta_{k+1}}{m_k} {\displaystyle \sum_{\ell=1}^{m_k}}c(\bmxi^\ell)
%\end{align*}
  \end{algorithmic}
\end{algorithm}

Similarly to Algorithm \ref{alg:rsag}, we can determine the convergence rate of Algorithm \ref{alg:rsag_specialized}.

\setcounter{theorem}{14}
\begin{theorem} \label{thm:convergence_multagent}
Suppose that Assumptions \ref{assumption: pvt_v2},  \ref{asp:compact_set}, \ref{asp:revenue_average}, and \ref{asp:choice_model} hold.
Moreover, suppose that Assumption~\ref{asp:continuous_dist} holds if $\bmxi$ is a continuous random vector.
%Let $L_{Ef}:= L_f + f_{\max} M$  and $m_k:=\left\lceil \frac{(c_{\max}M)^2k}{L_{Ef}\tilde{D}^2}\right\rceil$, where $c_{\max}=\max_{\bmxi\in \Xi} |c(\bmxi)|$. 
Let the inputs of Algorithm \ref{alg:rsag_specialized} other than $m_k$ be as in Theorem~\ref{thm:convergence} and let $m_k:=\left\lceil \frac{4(c_{\max}M)^2k}{L_{Ef}\tilde{D}^2}\right\rceil$, where $L_{Ef}:= L_f + f_{\max} M$.
Then, Algorithm \ref{alg:rsag_specialized} achieves the same convergence rate as in Theorem \ref{thm:convergence}.
\end{theorem}

\section{Experiments}
We conducted experiments on an application of \textbf{multiproduct pricing} to show that Algorithm \ref{alg:rsag_specialized} outputs solutions with higher total revenues compared with the existing methods. 
We performed synthetic experiments and simulation experiments with real retail data from a supermarket service provider in Japan.\footnote{We used publicly available data, ``New Product Sales Ranking'', provided by KSP-SP Co., Ltd, http://www.ksp-sp.com.} 
The details of our experiments are in Appenndix \ref{app:detail_exp}.

We implemented the following methods.
\\
\textbf{Proposed Method:} We implemented Algorithm \ref{alg:rsag_specialized} with $\alpha_k:=\frac{10}{k+1}, \beta_k:=\frac{0.1}{2m},\lambda_k:=\frac{k\beta_k}{2}$, $m_k=0.1km$, and $\zeta_k=\frac{1}{k}$, where $k$ is the current iteration number and $m$ is the number of buyers.
\\
\textbf{Proposed Method (fixed $\delta$):} This is the proposed method with a fixed $\delta_k$ from the information at the initial iterate $\bx_0$.
Specifically, $\delta_k$ was set to $\frac{1}{10^3}\sum_{\ell=1}^{10^3} c(\bmxi^\ell(\bx_0))$ for all $k$, where $\bmxi^\ell(\bx_0)\sim D(\bx_0)$.
\\
\textbf{Proposed Method ($\delta=0$):}
This is the proposed method with $\delta_k$ set to zero.
\\
\textbf{L2-Regularized Repeated Gradient Descent (L2-RGD($\alpha$)) \citep[Appendix E]{perdomo2020performative}:} 
\\
This method applies a repeated gradient descent \citep[Section 3.3]{perdomo2020performative} to the objective function with a regularization term $\frac{\alpha}{2}\|\bx-\bx_0\|$, where $\bx_0$ is the initial point.\footnote{This is introduced in \citep[Appendix E]{perdomo2020performative} as a remedy for the retraining method for non-strongly convex objective functions. 
Note that the retraining method was originally intended for strongly convex objective functions.}
We implemented this method for several $\alpha$.
\\
\textbf{Bayesian Optimization (BO) \citep{gpyopt2016}:}
This method sequentially searches for points where the objective value is likely to be small and outputs the solution with the lowest objective value among the evaluated points.
We used GPyOpt, a Python open-source library for Bayesian optimization \citep{gpyopt2016}.
\\
\textbf{Simultaneous Perturbation Stochastic Approximation (SPSA) \citep{spall2005introduction,spall1998implementation}:}
This method updates the current iterate by using approximated gradient, which calculated by the difference between objective values of two perturbed iterates.
%We set $c_k:=\frac{1}{(k+1)^{0.101}}$ and let $a_k:=\frac{0.16}{(100+k+1)^{0.602}}$ as the stepsize at each iteration \citep{spall1998implementation}
%We sample $\Delta^k$ from the Rademacher distribution, i.e. Bernoulli $\pm 1$ with probability $0.5$.
%The algorithm terminates at the iteration after $5m$ seconds or after $1000$ iterations.
\\
\textbf{Projected Sub-gradient Descent for Average Demand (PSD-AD) \citep[Section 3]{boyd2003subgradient}:} 
This is a projected subgradient descent method for deterministic pricing problems with average demand.
%We set the step size at each iteration so that the objective value decreases by repeatedly multiplying by $0.5$.
%The algorithm terminates at the iteration after $5m$ seconds or after $200$ iterations.

We performed our experiments under the following settings.\\
\textbf{Initial points.}
For all methods other than BO, we set the initial points as $\bx_0:=0.5\bme$, where $\bme \in \R^n$ is a vector with all elements equal to $1$.
BO first evaluates five random points; then it runs the Bayesian optimization.\\
\textbf{Metric.}
We computed $\frac{1}{10^3}\sum_{\ell=1}^{10^3} \left( -s(\bx_k,\bmxi^\ell(\bx_k)) + c(\bmxi^\ell(\bx_k))\right)$ for each iterate $\bx_k$, where $\bmxi^\ell(\bx_k) \sim D(\bx_k)$, and defined the smallest value among the iteration points as the \emph{Negative Expected Revenue} (NER).\\
\textbf{Termination criteria.}
We terminated all methods at a maximum computation time of $500$ seconds.

\subsection{Synthetic Experiments}
\paragraph{\textbf{Synthetic Parameter Setup.}}
We performed experiments by varying each parameter from
the following default settings.
We set $n:=20$ and $m:=200$, which are the numbers of products and buyers, respectively. 
For each product, we let the minimum price $x_{\min}$ be $0.01$ and the maximum price be $x_{\max}:=10$.
For the parameters of the function $p_i$, we generated $\alpha_i$ for each $i$ from a uniform distribution of $[0.01, 1]$, and we let $\gamma_i:=\frac{2\pi}{\sqrt{6}\alpha_i}$ and let $a_0:=0.25n$.
For the parameters of the function $c_i$ for each $i$, we set $\eta^1_i:=2.0w_i$, $\eta^2_i:=w_i$, and $\eta^3_i:=3.0w_i$, where $w_i$ was generated from a uniform distribution of $[0.25 \alpha_i, 0.5 \alpha_i]$. 
We let $l_i:=\frac{0.5m}{n}$ and $u_i:=\frac{1.5m}{n}$. 
We then varied $m$ and $n$ under these default settings.

\begin{table*}[t]
  \centering
  \caption{
  Results of synthetic experiments for 20 randomly generated problem instances. 
  The \emph{NER} (\emph{SD}) column represents the average (standard deviation) of NER.
  The best value of the average NER for each experiment is in bold. 
  }
  \label{tab:syn}
    \scalebox{0.7}{
      \begin{tabular}{crrrrrrrrrrrrrrrrrrrr} \toprule
($n$, $m$)
&\multicolumn{2}{c}{Proposed}&
\multicolumn{2}{c}{
\begin{tabular}{c}
Proposed \\ (fixed $\delta$)\end{tabular}}&
\multicolumn{2}{c}{
\begin{tabular}{c}
Proposed \\ ($\delta=0$)\end{tabular}}&
\multicolumn{2}{c}{
\begin{tabular}{c}
L2-RGD \\ ($\alpha=0.1$)\end{tabular}
}&
\multicolumn{2}{c}{
\begin{tabular}{c}
L2-RGD \\ ($\alpha=1$)\end{tabular}
}&
\multicolumn{2}{c}{
\begin{tabular}{c}
L2-RGD \\ ($\alpha=10$)\end{tabular}
}&
\multicolumn{2}{c}{BO}&
\multicolumn{2}{c}{SPSA}&
\multicolumn{2}{c}{PSD-AD}
\\ 
\cmidrule(lr){2-3} \cmidrule(lr){4-5} \cmidrule(lr){6-7} \cmidrule(lr){8-9} \cmidrule(lr){10-11} \cmidrule(lr){12-13} \cmidrule(lr){14-15}\cmidrule(lr){16-17} \cmidrule(lr){18-19}
&\multicolumn{1}{c}{NER}&\multicolumn{1}{c}{SD}& 
\multicolumn{1}{c}{NER}&\multicolumn{1}{c}{SD}&
\multicolumn{1}{c}{NER}&\multicolumn{1}{c}{SD}&
\multicolumn{1}{c}{NER}&\multicolumn{1}{c}{SD}&
\multicolumn{1}{c}{NER}&\multicolumn{1}{c}{SD}&
\multicolumn{1}{c}{NER}&\multicolumn{1}{c}{SD}&
\multicolumn{1}{c}{NER}&\multicolumn{1}{c}{SD}&
\multicolumn{1}{c}{NER}&\multicolumn{1}{c}{SD}&
\multicolumn{1}{c}{NER}&\multicolumn{1}{c}{SD} \\ 
\midrule
$(20, 200)$
&\textbf{-56.3}&5.4
&-54.9&5.5
&-54.7&5.2
&-28.0&10.6
&-28.2&10.7
&-28.2&10.7
&-22.3&6.2
&-33.7&17.6
&-45.9&8.9
\\
$(10, 200)$
&\textbf{-55.4}&8.4
&-54.5&8.8
&-54.4&8.3
&-8.2&23.0
&-9.8&23.2
&-11.9&23.4
&-34.4&10.6
&-46.8&13.4
&-31.8&13.9
\\ 
$(40, 200)$
&\textbf{-56.6}&3.2
&-54.7&3.5
&-54.1&3.4
&-24.5&9.5
&-24.4&9.5
&-24.4&9.6
&-14.6&3.9
&-1.7&13.0
&-47.4&3.7
\\
$(20, 100)$
&\textbf{-26.9}	&3.7
&-26.2	&3.8
&-26.1	&3.7
&-8.8	&6.8
&-8.8	&6.8
&-8.8	&6.8
&-11.7	&4.0
&-4.6	&6.3
&-20.4	&5.1
\\ 
$(20, 400)$
&\textbf{-106.9}	&14.2
&-104.4	&14.6
&-103.4	&14.0
&-38.5	&28.5
&-38.5	&28.2
&-39.2	&26.8
&-37.8	&7.9
&-36.4	&17.5
&-79.4	&21.0
\\ \bottomrule
\end{tabular}}
\end{table*}

\paragraph{Experimental Results} Table \ref{tab:syn} shows the results of the simulation experiments with different parameter values. The proposed method outperformed the baselines in terms of NER for all parameters, for the following reasons: (i) \textbf{Proposed (fixed $\delta$)} and \textbf{Proposed ($\delta=0$)} converged to low-quality local solutions because the variance of the gradient was larger than that of the proposed method; (ii) \textbf{L2-RGD} continued to increase prices without considering the effect of prices on the probability distribution, as shown in Example \ref{exa:exisiting_for_pricing} in Section \ref{subsec:related_decision_dependent_noise}, which led to unreasonably high prices; (iii) \textbf{BO} did not adequately explore $\bx$ because it took a lot of time to evaluate the objective value at each search point; (iv) \textbf{SPSA} did not accurately estimate the gradient because the noise in the gradient was too large; (v) \textbf{PSD-AD} ignored demand uncertainty, which increases the objective value since over/under demand occurs stochastically and causes unprofitable costs.

\subsection{Simulation Experiments with Real Data}

\paragraph{Data Set and Parameter Setup}
We used retail data from a supermarket service provider in Japan.
This data records the average sales prices of top-selling new products in food supermarkets.
We targeted sales data for $n=50$ different confectionery products for randomly selected weeks from 2022.
We set the recorded average selling price as the general value $\alpha_i$ for each product $i$.
The other parameters were set the same as in the synthetic experiment.
Since the parameter $w_i$ for each $i$ was generated randomly, experiments were performed on 20 problem instances for each week's data.

\begin{table*}[t]
  \centering
  \caption{
  Results of simulation experiments with real data for 20 randomly generated problem instances. 
The \emph{NER} (\emph{SD}) column represents the average (standard deviation) of NER.
The best value of the average NER for each experiment is in bold.
  }
  \label{tab:real}
    \scalebox{0.7}{
      \begin{tabular}{crrrrrrrrrrrrrrrrrrrr} \toprule
date
&\multicolumn{2}{c}{Proposed}&
\multicolumn{2}{c}{
\begin{tabular}{c}
Proposed \\ (fixed $\delta$)\end{tabular}}&
\multicolumn{2}{c}{
\begin{tabular}{c}
Proposed \\ ($\delta=0$)\end{tabular}}&
\multicolumn{2}{c}{
\begin{tabular}{c}
L2-RGD \\ ($\alpha=0.1$)\end{tabular}
}&
\multicolumn{2}{c}{
\begin{tabular}{c}
L2-RGD \\ ($\alpha=1$)\end{tabular}
}&
\multicolumn{2}{c}{
\begin{tabular}{c}
L2-RGD \\ ($\alpha=10$)\end{tabular}
}&
\multicolumn{2}{c}{BO}&
\multicolumn{2}{c}{SPSA}&
\multicolumn{2}{c}{PSD-AD}
\\ 
\cmidrule(lr){2-3} \cmidrule(lr){4-5} \cmidrule(lr){6-7} \cmidrule(lr){8-9} \cmidrule(lr){10-11} \cmidrule(lr){12-13} \cmidrule(lr){14-15}\cmidrule(lr){16-17} \cmidrule(lr){18-19}
&\multicolumn{1}{c}{NER}&\multicolumn{1}{c}{SD}& 
\multicolumn{1}{c}{NER}&\multicolumn{1}{c}{SD}&
\multicolumn{1}{c}{NER}&\multicolumn{1}{c}{SD}&
\multicolumn{1}{c}{NER}&\multicolumn{1}{c}{SD}&
\multicolumn{1}{c}{NER}&\multicolumn{1}{c}{SD}&
\multicolumn{1}{c}{NER}&\multicolumn{1}{c}{SD}&
\multicolumn{1}{c}{NER}&\multicolumn{1}{c}{SD}&
\multicolumn{1}{c}{NER}&\multicolumn{1}{c}{SD}&
\multicolumn{1}{c}{NER}&\multicolumn{1}{c}{SD} \\ 
\midrule
02/21--02/27
&\textbf{-28.1} 	&1.0 
&-21.5 	&1.4 
&-25.3 	&1.0 
&13.8 	&18.2 
&8.8 	&19.0 
&2.2 	&16.3 
&-8.5 	&2.1 
&10.3 	&8.7 
&-9.0 	&2.5 
\\
03/21--03/27
&\textbf{-20.6} 	&0.7 
&-20.1 	&1.0 
&-18.5 	&1.0 
&-7.5 	&3.4 
&-7.5 	&3.4 
&-7.6 	&3.4 
&-4.4 	&0.7 
&-10.4 	&3.3 
&-17.7 	&1.8 
\\ 
05/23--05/29
&\textbf{-22.6} 	&0.9 
&-17.8 	&1.8 
&-20.2 	&1.0 
&12.4 	&7.1 
&12.3 	&7.1 
&12.1 	&7.9 
&-6.1 	&1.6 
&-1.2 	&6.6 
&-10.2 	&3.2 
\\
06/20--06/26
&\textbf{-32.3} 	&2.1 
&-21.6 	&3.6 
&-28.8 	&2.4 
&79.2 	&39.5 
&55.0 	&58.9 
&53.3 	&57.9 
&-14.1 	&4.1 
&30.6 	&15.4 
&-8.1 	&5.6 
\\
08/08--08/14
&\textbf{-33.6} 	&0.9 
&-31.7 	&1.0 
&-31.2 	&1.0 
&-24.5 	&3.8 
&-24.5 	&3.7 
&-24.6 	&3.7 
&-7.2 	&1.8 
&-11.1 	&3.4 
&-29.4 	&1.6 
\\
09/19--09/25
&\textbf{-31.3} 	&1.5 
&-23.9 	&3.4 
&-28.5 	&1.8 
&0.0 	&24.4 
&-6.4 	&22.1 
&-10.3 	&18.5 
&-9.8 	&2.3 
&11.6 	&7.4 
&-13.5 	&5.3 
\\
12/05--12/11
&\textbf{-73.0} 	&3.2 
&-66.0 	&3.9 
&-71.1 	&3.4 
&172.3 	&30.7 
&152.6 	&44.0 
&146.2 	&37.2 
&-37.9 	&7.0 
&72.5 	&22.3 
&-28.9 	&10.8 
\\ \bottomrule
\end{tabular}}
\end{table*}

\paragraph{Experimental Results}
Table \ref{tab:real} shows the results of the experiments using real data from different weeks.
The proposed method was superior to the baseline in terms of NER for all weeks of data.

\section{Conclusion}
We formulated a new price optimization problem with decision-dependent uncertainty to address the drawbacks of existing formulations that (i) cannot deal with decision-dependent demand uncertainty, (ii) require discontinuous functions to define buyers' discrete actions, or (iii) have limited applications due to specific assumptions.
Moreover, we developed a projected stochastic gradient descent method by deriving an unbiased stochastic gradient with a variance reduction parameter.
Our method is guaranteed to converge to an $\epsilon$-stationary point.
Synthetic experiments and simulation experiments with real data confirmed the effectiveness of our formulation and method.

Our formulation and results suggest directions for further research.
The first is to construct methods to find a globally optimal solution rather than a stationary point (e.g., incorporating multi-start techniques \citep{Gyorgy2014multistart} into our methods or building fast Bayesian optimization under more specific assumptions).
The second is analyzing the performance of our method when some of our assumptions are relaxed.
This would include analyzing the performance when the probability density function is not differentiable and smoothed with the existing techniques.

\bibliographystyle{abbrvnat}
\bibliography{reference}

\appendix

\section{Proofs} \label{sec:proof}

\subsection{Proof of Proposition 1}

\begin{proof}
Assumption \ref{asp:compact_set} holds since $\mathcal C= [x_{\min},x_{\max}]^n$ and $\bmxi \in \{ 0,1,\dots,m\}^{n+1}$.
Therefore, we give proof for each of (i)--(iii) in Assumption~1. \\
{\bf (i)} From definitions of $s$ and $c$, the function $f(\bx, \bmxi)$ is differentiable w.r.t. $\bx$ and continuous w.r.t. $\bmxi$ for all 
$\bx \in \R^n$ and $\bmxi \in \Xi$. 
Moreover,
%\begin{align*}
    $\| \nabla_{\bx}f(\bx,\bmxi)\|=\| \nabla_{\bx} \left(-s(\bx,\bmxi)+c(\bmxi)\right)\| =\| -\nabla_{\bx}s(\bx,\bmxi)\| \le \sum_{i=1}^n |\xi_i| \le m$,
%\end{align*}
where the second inequality is due to the fact that the total demand for all products never exceeds the number $m$ of buyers.
Therefore, $f(\bx,\bmxi)$ is Lipschitz continuous with modulus $L_f=m$.

\noindent {\bf (ii)} Since $\Pr(\bmxi \mid \bx) = \prod_{i=0}^n {}_{m}C_{\xi_i} p_i(\bx)^{\xi_i}$, the function $\mpr(\bmxi \mid \bx)$ is differentiable w.r.t. $\bx$ and $\mpr(\bmxi \mid \bx) \neq 0$ for all $\bx \in \mR^n$ and $\bmxi \in \Xi$ from the definition of $p_i(\bx)$ for each $i\in \{0,1,\dots, n\}$.

\noindent {\bf (iii)}
We have $0 < p_i(\bx) < 1$ for all $\bx \in \mR^n$ and $i \in \{0,1,\dots,n\}$ from the definition of $p_i(\bx)$ for each $i\in \{0,1,\dots, n\}$.
Then, since 
%\begin{align*}
$\left(\nabla_{\bx} \mpr(\bmxi \mid \bx)\right)_k = \prod_{i=0}^n {}_{m}C_{\xi_i} p_i(\bx)^{\xi_i}\sum_{j=0}^n \frac{\xi_j}{p_j(\bx)} \frac{\partial p_j(\bx)}{\partial x_k}$,
%\end{align*}
we have 
\begin{align}
     \left(\frac{\nabla_{\bx} \mpr(\bmxi \mid \bx)}{\mpr(\bmxi \mid \bx)} \right)_k = \sum_{j=0}^n \frac{\xi_j}{p_j(\bx)} \frac{\partial p_j(\bx)}{\partial x_k}. \label{eq:partial_npr_pr}
\end{align}
Let $a_i(x):=e^{\gamma_i(\alpha_i-x)}$ for all $i \in I\ (=\{1,2,\dots,n\})$.
Then, $\frac{\partial a_i(x)}{\partial x}=-\gamma_i a_i(x)$. 
For $k \in I$, 
\small
\begin{align}
    \frac{\partial p_k(\bx)}{\partial x_k} &= \frac{\frac{\partial a_k(x_k)}{\partial x_k} (a_0+\sum_{i\in I} e^{\gamma_i(\alpha_i-x_i)}) -a_k(x_k) \frac{\partial a_k(x_k)}{\partial x_k}}{(a_0+\sum_{i\in I} e^{\gamma_i(\alpha_i-x_i)})^2} = \frac{-\gamma_k a_k(x_k) (a_0+\sum_{i\in I} e^{\gamma_i(\alpha_i-x_i)} -a_k(x_k))}{(a_0+\sum_{i \in I} e^{\gamma_i(\alpha_i-x_i)})^2} \nonumber\\
    &= -\gamma_k \frac{a_k(x_k)}{a_0+\sum_{i\in I} e^{\gamma_i(\alpha_i-x_i)}} \frac{a_0+\sum_{i\in I} e^{\gamma_i(\alpha_i-x_i)} -a_k(x_k)}{a_0+\sum_{i \in I} e^{\gamma_i(\alpha_i-x_i)}} =-\gamma_k p_k(\bx) (1-p_k(\bx)). \label{eq:partial_pi_xi}
\end{align}
\normalsize
For $k \in I$ and $j \in I \setminus \{k\}$,
\begin{align}
    \frac{\partial p_j(\bx)}{\partial x_k} &= \frac{-a_j(x_j)\frac{\partial a_k(x_k)}{\partial x_k}}{(a_0+\sum_{i\in I} e^{\gamma_i(\alpha_i-x_i)})^2} =\frac{-a_j(x_j)(-\gamma_k a_k(x_k))}{(a_0+\sum_{i \in I} e^{\gamma_i(\alpha_i-x_i)})^2}\nonumber\\
    &= \gamma_k\frac{a_j(x_j)}{a_0+\sum_{i \in I} e^{\gamma_i(\alpha_i-x_i)}} \frac{a_k(x_k)}{a_0+\sum_{i \in I} e^{\gamma_i(\alpha_i-x_i)}}=\gamma_k p_j(\bx) p_k(\bx).\label{eq:partial_pi_xj}
\end{align}
For $k \in I$,
\begin{align}
    \frac{\partial p_0(\bx)}{\partial x_k} &= \frac{-a_0\frac{\partial a_k(x_k)}{\partial x_k}}{(a_0+\sum_{i\in I} e^{\gamma_i(\alpha_i-x_i)})^2}=\frac{-a_0(-\gamma_k a_k(x_k))}{(a_0+\sum_{i \in I} e^{\gamma_i(\alpha_i-x_i)})^2}\nonumber\\
    &= \gamma_k\frac{a_0}{a_0+\sum_{i \in I} e^{\gamma_i(\alpha_i-x_i)}} \frac{a_k(x_k)}{a_0+\sum_{i \in I} e^{\gamma_i(\alpha_i-x_i)}}=\gamma_k p_0(\bx) p_k(\bx).\label{eq:partial_p0_xj}
\end{align}
From \eqref{eq:partial_npr_pr}--\eqref{eq:partial_p0_xj},
\small
\begin{align*}
     \left(\frac{\nabla_{\bx} \mpr(\bmxi \mid \bx)}{\mpr(\bmxi \mid \bx)} \right)_k 
     &= \sum_{j\in I \setminus \{k\}} \left( \frac{\xi_j}{p_j(\bx)} \gamma_k p_k(\bx) p_j(\bx)\right) -\frac{\xi_k}{p_k(\bx)} \gamma_k p_k(\bx) (1-p_k(\bx))+\frac{\xi_0}{p_0(\bx)} \gamma_k p_0(\bx) p_k(\bx) \\
     &=\sum_{j\in I \cup \{0\} \setminus \{k\}} \xi_j\gamma_k p_k(\bx) -\xi_k\gamma_k(1-p_k(\bx))\le \sum_{j \in I\cup \{0\}} | \xi_j | |\gamma_k| \le m\gamma^{\max},
\end{align*}
\normalsize
where the first inequality follows from $0 < p_k(\bx) < 1$ for all $\bx \in \mR^n$. 
The second inequality follows from the definition whereby $\sum_{j \in I\cup \{0\}} | \xi_j |$ is equal to the number $m$ of buyers.
Then, for all $x \in \mR^n$ and $\bmxi \in \Xi$,
%\begin{align*}
    $\left\| \frac{\nabla_{\bx} \mpr(\bmxi \mid \bx)}{\mpr(\bmxi \mid \bx)} \right\| \le \sum_{k=1}^n \left|\left(\frac{\nabla_{\bx} \mpr(\bmxi \mid \bx)}{\mpr(\bmxi \mid \bx)} \right)_k \right| \le n m \gamma^{\max}$.  
%\end{align*}
\end{proof}

\subsection{Proof of Proposition 2}

\begin{proof}
Assumption \ref{asp:compact_set} holds since $\mathcal C= [x_{\min},x_{\max}]^I$ and $\xi_i \in \{ 0,1,\dots, d_i\}$ for all $i \in I$.
Therefore, we give proof for each of (i)--(iii) in Assumption~1. \\
{\bf (i)} From the definition, the value of $f(\bx,\bmxi)$ is independent of $\bx$. 
Therefore, $f(\bx,\bmxi)$ is differentiable and Lipschitz continuous with modulus $L_f=0$ w.r.t. $\bx$.
Moreover, $f(\bx,\bmxi)$ is continuous w.r.t. $\bmxi$ from the definition since $q_H$, $q_R$, and $k$ are continuous functions.

\noindent {\bf (ii)} Since $\Pr(\bmxi \mid \bx) =\prod_{i \in I} {}_{d_i}C_{\xi_i} p_i(x_i)^{\xi_i}(1-p_i(x_i))^{d_i-\xi_i}$, $\mpr(\bmxi \mid \bx)$ is differentiable w.r.t. $\bx$ from the definition of $p_i$.
Moreover, since $0 < p_i(x) < 1$ for all $x \in \mR$ and $i \in I$, $\mpr(\bmxi \mid \bx) \neq 0$ for all $(\bx, \bmxi) \in \mR^n \times \Xi$.

\noindent {\bf (iii)}
For all $x \in \mR$ and $i \in I$, 
the definition of $p_i$ gives $0 < p_i(x) < 1$.
Then, since 
\begin{align*}
     \left(\nabla_{\bx} \mpr(\bmxi \mid \bx)\right)_k & = \prod_{i\in I} {}_{d_i}C_{\xi_i} p_i(x_i)^{\xi_i} (1-p_i(x_i))^{d_i-\xi_i}\left( \frac{\xi_k p'_k(x_k)}{p_k(x_k)} - \frac{(d_k-\xi_k) p'_k(x_k)}{1-p_k(x_k)} \right),
\end{align*}
we have 
\begin{align}
     \left(\frac{\nabla_{\bx} \mpr(\bmxi \mid \bx)}{\mpr(\bmxi \mid \bx)} \right)_k =  \frac{\xi_k p'_k(x_k)}{p_k(x_k)} - \frac{(d_k-\xi_k) p'_k(x_k)}{1-p_k(x_k)} . \label{eq:partial_npr_pr_2}
\end{align}    
Here,
\begin{align}
    p'_k(x)&= \frac{-\beta_k e^{\alpha_k h_k+\beta_kx+\gamma_k}}{(1+e^{\alpha_k h_k+\beta_kx+\gamma_k})^2} =-\beta_k \frac{1}{(1+e^{\alpha_k h_k+\beta_kx+\gamma_k})} \frac{e^{\alpha_k h_k+\beta_kx+\gamma_k}}{(1+e^{\alpha_k h_k+\beta_kx+\gamma_k})} =-\beta_k p_k(x) (1-p_k(x)). \label{eq:partial_pi_xi_2}
\end{align}
From \eqref{eq:partial_npr_pr_2} and \eqref{eq:partial_pi_xi_2},
\begin{align*}
     \left(\frac{\nabla_{\bx} \mpr(\bmxi \mid \bx)}{\mpr(\bmxi \mid \bx)} \right)_k &=-\beta_k (1-p_k(x_k))\xi_k +\beta_k p_k(x_k)(d_k-\xi_k)=\beta_k (d_kp_k(x_k)-\xi_k)\\
     &\le |\beta_k| |d_kp_k(x_k)- \xi_k|\le |\beta_k| d_k,
\end{align*}
where the second inequality follows from $\xi_k \in \{0,1,\dots,d_k\}$ and $0 \le p_k(x_k) \le 1$ for all $k \in I$.
Then, for all $x \in \mR^I$ and $\bmxi \in \Xi$,
%\begin{align*}
    $\left\| \frac{\nabla_{\bx} \mpr(\bmxi \mid \bx)}{\mpr(\bmxi \mid \bx)} \right\| \le  \sum_{k\in I} \left|\left(\frac{\nabla_{\bx} \mpr(\bmxi \mid \bx)}{\mpr(\bmxi \mid \bx)} \right)_k \right|\le |I| \max_{i \in I}\left( |\beta_i| d_i \right)$. 
%\end{align*}
\end{proof}

\subsection{Proof of Proposition 3}
\begin{proof}
Assumption \ref{asp:compact_set} holds since $\mathcal C= [x_{\min},x_{\max}]^n$ and $\Xi = \{\bmxi \mid \forall i \in [n], 0 \le \xi_i \le \xi_i^{\max} \}$.
Moreover, Assumption \ref{asp:continuous_dist} holds since
$\Xi$ is a Borel set on $\mR^n$ and $\mpr(\bmxi\mid\bx)$ is continuous w.r.t. $\bmxi$ for all $\bx \in \R^n$ from the definition of $\mpr(\bmxi\mid\bx)$.
We give proof for each condition of Assumption~1. \\
{\bf (i)} From definitions of $s$ and $c$,  $f(\bx, \bmxi)$ is differentiable w.r.t. $\bx$ and continuous w.r.t. $\bmxi$ for all $\bx \in \R^n$ and $\bmxi \in \Xi$. Moreover, $\| \nabla_{\bx}f(\bx,\bmxi)\|=\| \bmxi\| \le n\xi^{\max}$.
Therefore, $f(\bx,\bmxi)$ is Lipschitz continuous with modulus $L_f=n\xi^{\max}$.

\noindent {\bf (ii)} $\mpr(\bmxi \mid \bx)$ is differentiable w.r.t. $\bx$ and $\mpr(\bmxi \mid \bx) \neq 0$ for all $\bx \in \mR^n$ and $\bmxi \in \Xi$ from definitions of $\Pr(\bmxi \mid \bx)$, $C^i(\bx)$, and $\bv^i(\bx)$.

\noindent {\bf (iii)}
Let $g_i(\bx,\xi):=-\frac{(\xi-\bv^i(\bx)^\top \ba^i)^2}{2((\sigma^i)^2-\bv^i(\bx)^\top A^i \bv^i(\bx))}$ and $h_i(\bx):=(\sigma^i)^2-\bv^i(\bx)^\top A^i \bv^i(\bx)$.
Then, 
$\Pr(\bmxi \mid \bx):= 
\prod_{i=1}^n \frac{1}{C^i(\bx)\sqrt{2\pi h_i(\bx)}} \exp\left(g_i(\bx,\xi_i)\right)$.
Therefore,
\small
\begin{align}
\frac{\partial \Pr(\bmxi \mid \bx)}{\partial x_k}=\prod_{i=1}^n \frac{1}{C^i(\bx)\sqrt{2\pi h_i(\bx)}} \exp\left(g_i(\bx,\xi_i)\right) \left(-\sum_{i=1}^n \frac{1}{C^i(\bx)} \frac{\partial C^i(\bx)}{\partial x_k} -\sum_{i=1}^n\frac{1}{2}\frac{1}{h^i(\bx)} \frac{\partial h^i(\bx)}{\partial x_k}  + \sum_{i=1}^n \frac{\partial g_i(\bx,\xi_i)}{\partial x_k} \right). \label{eq:pr_partial_xk}
\end{align}
\normalsize
Here,
\begin{align}
\left|\frac{\partial v_j^i(\bx)}{\partial x_k}\right| =\left|-\theta_1^i\frac{2x_k-2\hat{x}_k^j}{\theta_2^i} \exp \left(-\frac{\|\bx-\hat{\bx}^j\|^2}{\theta_2^i}\right)\right| \le \frac{2|\theta_1^i||x_k-\hat{x}_k^j|}{|\theta_2^i|} \le 
\frac{2\theta_1^{\max}(x_{\max}-x_{\min})}{\theta_2^{\min}}. \label{eq:partial_v_xk_bound} 
\end{align}
For all $\xi_i \in \R$,
\small
\begin{align}
\left| \frac{\partial g_i(\bx,\xi_i)}{\partial x_k} \right|&=
\left|-\frac{2(\xi_i-\bv^i(\bx)^\top \ba^i)(-\sum_{j=1}^N \frac{\partial v_j^i(\bx)}{\partial x_k}a_j^i)}{2((\sigma^i)^2-\bv^i(\bx)^\top A^i \bv^i(\bx))} \right. \nonumber \\
& \quad \left. + \frac{(\xi_i-\bv^i(\bx)^\top \ba^i)^2(-2\sum_{s=1}^N\sum_{t=1}^N A_{st}^i (v_t^i(\bx)\frac{\partial v_s^i(\bx)}{\partial x_k}+v_s^i(\bx)\frac{\partial v_t^i(\bx)}{\partial x_k}))}{(2((\sigma^i)^2-\bv^i(\bx)^\top A^i \bv^i(\bx)))^2}\right| \nonumber\\
&\le \left|\frac{\xi_i-\bv^i(\bx)^\top \ba^i}{(\sigma^i)^2-\bv^i(\bx)^\top A^i \bv^i(\bx)}\right| \left|\sum_{j=1}^N\frac{\partial v_j^i(\bx)}{\partial x_k}a_j^i\right| \nonumber\\
& \quad + \left|\frac{(\xi_i-\bv^i(\bx)^\top \ba^i)^2}{(2((\sigma^i)^2-\bv^i(\bx)^\top A^i \bv^i(\bx)))^2}\right| \left|2\sum_{s=1}^N\sum_{t=1}^N A_{st}^i \left(v_t^i(\bx)\frac{\partial v_s^i(\bx)}{\partial x_k}+v_s^i(\bx)\frac{\partial v_t^i(\bx)}{\partial x_k}\right)\right| \nonumber\\
& \le \frac{\xi^{\max} + N\theta_1^{\max} a^{\max}}{\Delta} \cdot 2Na^{\max}\frac{\theta_1^{\max}(x_{\max}-x_{\min})}{\theta_2^{\min}} \nonumber\\
& \quad + \frac{(\xi^{\max} + N\theta_1^{\max} a^{\max})^2}{4\Delta^2}\cdot 8 N^2 A^{\max}\frac{(\theta_1^{\max})^2 (x_{\max}-x_{\min})}{\theta_2^{\min}} \nonumber\\
& \le \frac{2N\theta_1^{\max} (\xi^{\max} + N\theta_1^{\max} a^{\max})(x_{\max}-x_{\min})}{\Delta \theta_2^{\min}}\left(a^{\max}+NA^{\max}\theta_1^{\max}\frac{\xi^{\max} + N\theta_1^{\max} a^{\max}}{\Delta}\right), \label{eq:partial_g_xk_bound}
\end{align}
\normalsize
where the second inequality comes from  \eqref{eq:partial_v_xk_bound} and the fact that $(\sigma^i)^2-\bv^i(\bx)^\top A^i \bv^i(\bx) \ge \Delta$ and $|v_t^i(\bx)|\le |\theta_1^i| \le \theta_1^{\max}$ for $i \in [n]$ and $t\in [N]$.

Moreover, for all $\bx \in \R^n$,
\begin{align}
\left|\frac{1}{2} \frac{1}{h_i(\bx)}\frac{\partial h_i(\bx)}{\partial x_k} \right|&=\frac{1}{2}\left|\frac{1}{h_i(\bx)}\right| \left|\sum_{s=1}^N\sum_{t=1}^N A_{st}^i \left(v_t^i(\bx)\frac{\partial v_s^i(\bx)}{\partial x_k}+v_s^i(\bx)\frac{\partial v_t^i(\bx)}{\partial x_k}\right)\right| \nonumber\\
&\le \frac{2}{\Delta}N^2 A^{\max}\frac{(\theta_1^{\max})^2 (x_{\max}-x_{\min})}{\theta_2^{\min}}, \label{eq:partial_hi_bounded}
\end{align}
where the inequality follows from  \eqref{eq:partial_v_xk_bound} and the fact that $h_i(\bx) \ge \Delta$ and $|v_t^i(\bx)| \le \theta_1^{\max}$ for for $i \in [n]$ and $t\in [N]$.

Let $r_1:=\frac{2N\theta_1^{\max} (\xi^{\max} + N\theta_1^{\max} a^{\max})(x_{\max}-x_{\min})}{\Delta \theta_2^{\min}}\left(a^{\max}+NA^{\max}\theta_1^{\max}\frac{\xi^{\max} + N\theta_1^{\max} a^{\max}}{\Delta}\right)$,
and \\
$r_2:= \frac{2}{\Delta}N^2 A^{\max}\frac{(\theta_1^{\max})^2 (x_{\max}-x_{\min})}{\theta_2^{\min}}$. 
Then,
\begin{align}
\left|\frac{1}{C^i(\bx)} \frac{\partial C^i(\bx)}{\partial x_k}\right| &= \left|\frac{1}{C^i(\bx)} \frac{\partial \int_{0}^{\xi_i^{\max}} \frac{1}{\sqrt{2\pi h_i(\bx)}} \exp\left(g_i(\bx,\phi)\right) d\phi}{\partial x_k}\right| \nonumber\\
&= \left|\frac{1}{C^i(\bx)} \int_{0}^{\xi_i^{\max}} \frac{1}{\sqrt{2\pi h_i(\bx)}} \exp\left(g_i(\bx,\phi)\right)  \left(-\frac{1}{2}\frac{1}{h_i(\bx)}\frac{\partial h_i(\bx)}{\partial x_k} +\frac{\partial g_i(\bx,\phi)}{\partial x_k}\right) d\phi\right| \nonumber\\
&\le \left|\frac{1}{C^i(\bx)} \int_{0}^{\xi_i^{\max}} \frac{1}{\sqrt{2\pi h_i(\bx)}} \exp\left(g_i(\bx,\phi)\right) 
\left( \left| \frac{1}{2}\frac{1}{h_i(\bx)}\frac{\partial h_i(\bx)}{\partial x_k} \right| +\left|\frac{\partial g_i(\bx,\phi)}{\partial x_k}\right| \right) d\phi\right| \nonumber\\
&\le \left|\frac{1}{C^i(\bx)} \int_{0}^{\xi_i^{\max}} \frac{1}{\sqrt{2\pi h_i(\bx)}} \exp\left(g_i(\bx,\phi)\right) 
\left(r_1+r_2 \right) d\phi\right|  
=  r_1 + r_2, \nonumber
\end{align}
where the second inequality follows from \eqref{eq:partial_g_xk_bound} and \eqref{eq:partial_hi_bounded}.
Here, from \eqref{eq:pr_partial_xk},
\small
\begin{align*}
\left|\frac{\partial \Pr(\bmxi \mid \bx)}{\partial x_k} \cdot \frac{1}{\Pr(\bmxi \mid \bx)}\right| &\le \sum_{i=1}^n \left| \frac{1}{C^i(\bx)} \frac{\partial C^i(\bx)}{\partial x_k}\right| + \sum_{i=1}^n  \left| \frac{1}{2}\frac{1}{h^i(\bx)}\frac{\partial h^i(\bx)}{\partial x_k} \right| +\sum_{i=1}^n \left| \frac{\partial g_i(\bx,\xi_i)}{\partial x_k} \right| \le 2n(r_1 + r_2).
\end{align*}
\normalsize
Therefore,
%\begin{align*}
$\left\|\frac{\nabla \Pr(\bmxi \mid \bx)}{\Pr(\bmxi \mid \bx)} \right\| 
\le \sum_{k=1}^n\left|\frac{\partial \Pr(\bmxi \mid \bx)}{\partial x_k} \cdot \frac{1}{\Pr(\bmxi \mid \bx)}\right| 
\le 2n^2 (r_1 + r_2)$.
%\end{align*}
Condition (iii) of Assumption \ref{assumption: pvt_v2} holds from the definition of $r_1$, $r_2$, and $M$.
  \end{proof}

\subsection{Proof of Lemma \ref{lem:int_lim_change}}
\begin{proof}
For a given $\bx \in C$, let $\{\Delta_k\}$ be a sequence of scalars such that $\lim_{k \to \infty} \Delta_k = 0$ and $\bx+\Delta_k\bme^i \in \mathcal C$, where $\bme^i$ is a vector such that the $i$-th element is $1$ and other elements are $0$. 
Let $g_{k,i}(\bx,\bmxi):=\frac{h(\bx+\Delta_k\bme^i,\bmxi)-h(\bx,\bmxi)}{\Delta_k}$.
There exists $\bx' \in \mathcal C$ such that $g_{k,i}(\bx,\bmxi)=\frac{\partial h(\bx',\bmxi)}{\partial x_i}$ from the mean-value theorem.
Moreover, $f_{\max}(=\max_{\bx\in \mathcal C, \bmxi \in \Xi} |f(\bx,\bmxi)|)$ exists since $\Xi$ and $\mathcal C$ are compact from Assumption 2 and $f$ is a real-valued continuous function from Assumption 1.
Then, for all $\bx \in \mathcal \mathcal C$ and $i=1,\dots,n$,
\begin{align*}
|g_{k,i}(\bx,\bmxi)| &= \left|\frac{\partial h(\bx',\bmxi)}{\partial x_i}\right|= \left|\frac{\partial f(\bx',\bmxi)}{\partial x_i} \mpr(\bmxi\mid\bx')+f(\bx',\bmxi)
\frac{\partial \mpr(\bmxi \mid \bx')}{\partial x_i}\right| \\
&\le \left|\frac{\partial f(\bx',\bmxi)}{\partial x_i}\right|+f_{\max}\left|\frac{\partial \mpr(\bmxi \mid \bx')}{\partial x_i}\frac{1}{\mpr(\bmxi \mid \bx')}\right|\le L_f + f_{\max}M,
\end{align*}
where the first inequality comes from $0 < \mpr(\bmxi\mid\bx') \le 1$, and the second inequality follows from conditions (i) and (iii) of Assumption~1.
Here, $g_{k,i}$ is measurable on $\Xi$ since $\Xi$ is a Borel set and $g_{k,i}$ is continuous \emph{w.r.t.} $\bmxi$ from Assumption~3 and the definitions of $g_{k,i}$ and $h$.
The constant function $r(\bmxi):=L_f + f_{\max}M < \infty$ is integrable over $\Xi$.
Moreover, $g_{k,i}(\bx,\bmxi) \to \frac{\partial h(\bx,\bmxi)}{\partial x_i}$ pointwise when $k \to \infty$ since $h(\bx,\bmxi)$ is differentiable \emph{w.r.t.} $\bx$ from conditions (i) and (ii) of Assumption 1.
Then, the Lebesgue dominated convergence theorem \citep[Chapter 4.4, page 88]{royden1988real} holds for $g_{k,i}$ for all $\bx \in \mathcal{C}$ and $i=1,\dots, n$, that is,
\begin{align}
\lim_{k\to \infty} \int_{\bmxi \in \Xi} g_{k,i}(\bx,\bmxi) d\bmxi = \int_{\bmxi \in \Xi} \lim_{k\to \infty} g_{k,i}(\bx,\bmxi) d\bmxi,\ \textrm{for all} \  \bx \in \mathcal{C}\ \textrm{and}\  i=1,\dots, n. \nonumber
\end{align}
Then, for all $\bx \in \mathcal{C}$ and $i=1,\dots,n$,
\small
\begin{align}
\left(\nabla_{\bx} \int_{\bmxi \in \Xi} h(\bx,\bmxi) d\bmxi\right)_i \nonumber
&= \lim_{k \to \infty} 
\frac{\int_{\bmxi \in \Xi} h(\bx+\Delta_k\bme^i,\bmxi) d\bmxi - \int_{\bmxi \in \Xi} h(\bx,\bmxi) d\bmxi}{\Delta_k} \\
&=\lim_{k \to \infty} 
\int_{\bmxi \in \Xi} \frac{h(\bx+\Delta_k\bme^i,\bmxi)- h(\bx,\bmxi)}{\Delta_k} d\bmxi \nonumber\\
&=\lim_{k\to \infty} \int_{\bmxi \in \Xi} g_{k,i}(\bx,\bmxi) d\bmxi = \int_{\bmxi \in \Xi} \lim_{k\to \infty} g_{k,i}(\bx,\bmxi) d\bmxi \nonumber\\
%=\int_{\bmxi \in \Xi} \lim_{k \to \infty} \frac{h(\bx+\Delta_k\bme^i,\bmxi)- h(\bx,\bmxi)}{\Delta_k} d\bmxi \nonumber\\
&= \int_{\bmxi \in \Xi}\left(\nabla_{\bx}  h(\bx,\bmxi) \right)_id\bmxi.\nonumber
\end{align}
\normalsize
Therefore, for all $\bx \in \mathcal{C}$,
\begin{align*}
\nabla_{\bx} \int_{\bmxi \in \Xi} h(\bx,\bmxi) d\bmxi = \int_{\bmxi \in \Xi} \nabla_{\bx} h(\bx,\bmxi)d\bmxi.  \end{align*}
\end{proof}

\subsection{Proof of Lemma \ref{lem:int_lim_change_2}}
\begin{proof}
For given $\bx \in C$, let $\{\Delta_k\}$ be a sequence of scalars such that $\lim_{k \to \infty} \Delta_k = 0$ and $\bx+\Delta_k\bme^i \in \mathcal C$, where $\bme^i$ is a vector such that the $i$-th element is $1$ and other elements are $0$.
Let $g_{k,i}(\bx,\bmxi):=\frac{q(\bmxi)\mpr(\bmxi \mid \bx+\Delta_k\bme^i)-q(\bmxi)\mpr(\bmxi \mid \bx)}{\Delta_k}$.
There exists $\bx' \in \mathcal C$ such that $g_{k,i}(\bx,\bmxi)=\frac{\partial q(\bmxi)\mpr(\bmxi \mid \bx')}{\partial x_i}$ from the mean-value theorem.
Moreover, let $q^{\max}:=\max_{\bmxi \in \Xi}|q(\bmxi)|$, which exists since $\Xi$ is compact from Assumption 2 and $q$ is a real-valued continuous function.
Then, for all $\bx \in \mathcal \mathcal C$ and $i=1,\dots,n$,
\begin{align*}
|g_{k,i}(\bx,\bmxi)| &= \left|\frac{\partial q(\bmxi)\mpr(\bmxi \mid \bx')}{\partial x_i}\right| = |q(\bmxi)|\left|
\frac{\partial \mpr(\bmxi \mid \bx')}{\partial x_i}\right| \\
&\le |q(\bmxi)|\left|
\frac{\partial \mpr(\bmxi \mid \bx')}{\partial x_i} \frac{1}{\mpr(\bmxi \mid \bx')} \right| \le q^{\max} M,
\end{align*}
where the first inequality follows from $0 < \mpr(\bmxi\mid\bx') \le 1$.
The second inequality comes from condition (iii) of Assumption 1.
Here, $g_{k,i}$ is measurable on $\Xi$ since $\Xi$ is a Borel set and $g_{k,i}$ is continuous \emph{w.r.t.} $\bmxi$ from Assumption~3 and the definition of $g_{k,i}$.
The constant function $r(\bmxi):=q^{\max}M < \infty$ is integrable over $\Xi$.
Moreover, $g_{k,i}(\bx,\bmxi) \to \frac{\partial q(\bmxi)\mpr(\bmxi \mid \bx)}{\partial x_i}$ pointwise when $k \to \infty$ since $\mpr(\bmxi \mid \bx)$ is differentiable \emph{w.r.t.} $\bx$ from condition (ii) of Assumption~1.
Then, the Lebesgue dominated convergence theorem \citep[Chapter 4.4, page 88]{royden1988real} holds for $g_{k,i}$ for all $\bx \in \mathcal{C}$ and $i=1,\dots, n$, that is,
\begin{align}
\lim_{k\to \infty} \int_{\bmxi \in \Xi} g_{k,i}(\bx,\bmxi) d{\bmxi} = \int_{\bmxi \in \Xi} \lim_{k\to \infty} g_{k,i}(\bx,\bmxi) d{\bmxi},\ \textrm{for all} \  \bx \in \mathcal{C}\ \textrm{and}\  i=1,\dots, n. \nonumber
\end{align}
Then, for all $\bx \in \mathcal{C}$ and $i=1,\dots,n$,
\begin{align}
&\left(\nabla_{\bx} \int_{\bmxi \in \Xi} q(\bmxi)\mpr(\bmxi \mid \bx) d\bmxi\right)_i \nonumber\\
&= \lim_{k \to \infty} 
\frac{\int_{\bmxi \in \Xi} q(\bmxi)\mpr(\bmxi \mid \bx+\Delta_k\bme^i) d\bmxi - \int_{\bmxi \in \Xi} q(\bmxi)\mpr(\bmxi \mid \bx) d\bmxi}{\Delta_k}\nonumber\\
&=
\lim_{k \to \infty} 
\int_{\bmxi \in \Xi} \frac{q(\bmxi)\mpr(\bmxi \mid \bx+\Delta_k\bme^i)-q(\bmxi) \mpr(\bmxi \mid \bx)}{\Delta_k} d\bmxi \nonumber\\
&=\lim_{k\to \infty} \int_{\bmxi \in \Xi} g_{k,i}(\bx,\bmxi) d\bmxi = \int_{\bmxi \in \Xi} \lim_{k\to \infty} g_{k,i}(\bx,\bmxi) d\bmxi\nonumber\\
%&=\int_{\bmxi \in \Xi} \lim_{k \to \infty} \frac{q(\bmxi)\mpr(\bmxi \mid \bx+\Delta_k\bme^i)- q(\bmxi)\mpr(\bmxi \mid \bx)}{\Delta_k} d\bmxi \nonumber\\
&= \int_{\bmxi \in \Xi}\left(\nabla_{\bx} ( q(\bmxi)\mpr(\bmxi \mid \bx)) \right)_id\bmxi.
\end{align}
Therefore, for all $\bx \in \mathcal{C}$,
\begin{align*}
\nabla_{\bx} \int_{\bmxi \in \Xi} q(\bmxi)\mpr(\bmxi \mid \bx) d\bmxi = \int_{\bmxi \in \Xi} \nabla_{\bx} (q(\bmxi)\mpr(\bmxi \mid \bx))d{\bmxi}.\nonumber
\end{align*}
\end{proof}

\subsection{Proof of Lemma \ref{sfo_lemma}}
\begin{proof}
We have
\begin{align}
 &\int_{\bmxi \in \Xi} \left( \delta \frac{\nabla_{\bx} \mpr(\bmxi\mid\bx)}{\mpr(\bmxi\mid\bx)}\right) \mpr(\bmxi\mid\bx) d\bmxi =\delta \int_{\bmxi \in \Xi} \nabla_{\bx} \mpr(\bmxi\mid\bx) d\bmxi \nonumber \\
  &=\delta \nabla_{\bx} \int_{\bmxi \in \Xi} \mpr(\bmxi\mid\bx) d\bmxi=\delta \nabla_{\bx}(1)=0, \label{eq:delta_int_0}
\end{align}
where the second equality comes from Lemma \ref{lem:int_lim_change_2} with $q(\bmxi) = 1$.
Then,
\small
\begin{align*}
 &\nabla_{\bx} \E_{\bmxi \sim D(\bx)}[f(\bx,\bmxi)] = \nabla_{\bx} \int_{\bmxi \in \Xi} f(\bx,\bmxi)  d\mpr(\bmxi\mid\bx) = \nabla_{\bx} \int_{\bmxi \in \Xi} f(\bx,\bmxi) \mpr(\bmxi\mid\bx) d\bmxi \\
 &=\int_{\bmxi \in \Xi} \nabla_{\bx} \left( f(\bx,\bmxi) \mpr(\bmxi\mid\bx)\right) d\bmxi =\int_{\bmxi \in \Xi} \nabla_{\bx} f(\bx,\bmxi) \mpr(\bmxi\mid\bx) + f(\bx,\bmxi) \nabla_{\bx} \mpr(\bmxi\mid\bx)  d\bmxi \\
 &=\int_{\bmxi \in \Xi} \left(\nabla_{\bx} f(\bx,\bmxi) +  f(\bx,\bmxi) \frac{\nabla_{\bx} \mpr(\bmxi\mid\bx)}{\mpr(\bmxi\mid\bx)}\right) \mpr(\bmxi\mid\bx) d\bmxi \\
  &=\int_{\bmxi \in \Xi} \left(\nabla_{\bx} f(\bx,\bmxi) +  f(\bx,\bmxi) \frac{\nabla_{\bx} \mpr(\bmxi\mid\bx)}{\mpr(\bmxi\mid\bx)}\right) d\mpr(\bmxi\mid\bx) \\
    &=\int_{\bmxi \in \Xi} \left(\nabla_{\bx} f(\bx,\bmxi) +  (f(\bx,\bmxi)-\delta) \frac{\nabla_{\bx} \mpr(\bmxi\mid\bx)}{\mpr(\bmxi\mid\bx)}\right) d\mpr(\bmxi\mid\bx) \\
 &=\E_{\bmxi \sim D(\bx)} \left[\nabla_{\bx} f(\bx,\bmxi)+(f(\bx,\bmxi)-\delta) \frac{\nabla_{\bx} \mpr(\bmxi\mid\bx)}{\mpr(\bmxi\mid\bx)}\right].
\end{align*}
\normalsize
Here, the third equality obviously holds when $\bmxi$ is a discrete random vector.
If $\bmxi$ is a continuous random vector, the third equality follows from Lemma \ref{lem:int_lim_change} since Assumptions \ref{assumption: pvt_v2}--\ref{asp:continuous_dist} hold.
The fourth equality is due to the fact that $f(\bx,\bmxi)$ and $\mpr(\bmxi\mid\bx)$ are differentiable w.r.t. $\bx$ from conditions (i) and (ii) of Assumption~\ref{assumption: pvt_v2}.
The fifth equality is due to the fact that $\mpr(\bmxi\mid\bx) \neq 0$ from condition (ii) of Assumption~\ref{assumption: pvt_v2}.
The seventh equality comes from \eqref{eq:delta_int_0}.
Then, Lemma \ref{sfo_lemma} holds from Definition~2.
\end{proof}

\subsection{Proof of Lemma \ref{lemma:SFO_variance}}
\begin{proof}
We have 
\begin{align*}
&\E_{\bmxi' \sim D(\bx)}[\| g_2(\bx,\bmxi') - \nabla_{\bx} \E_{\bmxi \sim D(\bx)}[f(\bx,\bmxi)]  \|^2] 
= \E_{\bmxi' \sim D(\bx)}[\| g_2(\bx,\bmxi')\|^2] - \|\nabla_{\bx} \E_{\bmxi \sim D(\bx)}[f(\bx,\bmxi)]\|^2 \\
&\le \E_{\bmxi' \sim D(\bx)}[\| g_2(\bx,\bmxi')\|^2] \le \E_{\bmxi' \sim D(\bx)}\left[\left(\|\nabla_{\bx} f(\bx,\bmxi')\|+|f(\bx,\bmxi')-\delta|\left\|\frac{\nabla_{\bx} \mpr(\bmxi'\mid\bx)}{\mpr(\bmxi'\mid\bx)}\right\|\right)^2\right]  \\
&\le \E_{\bmxi' \sim D(\bx)}[(L_f+2f_{\max}M)^2] = (L_f+2f_{\max}M)^2,
\end{align*}
where the first equality comes from Lemma \ref{sfo_lemma}.
The third inequality follows 
from conditions (i) and (iii) of Assumption~1 and Assumption \ref{asp:compact_set}.
  \end{proof}

\subsection{Proof of Lemma \ref{lem:entire_derivative_bound}}
\begin{proof}
    For all $\bx \in \mathcal C$,
    \begin{align*}
    &\|\nabla_{\bx} \E_{\bmxi \sim D(\bx)}[f(\bx,\bmxi)] \| 
    = \left\|\E_{\bmxi \sim D(\bx)} \left[\nabla_{\bx} f(\bx,\bmxi)+f(\bx,\bmxi) \frac{\nabla_{\bx} \mpr(\bmxi\mid\bx)}{\mpr(\bmxi\mid\bx)}\right]\right\| \\
    & \le \left\|\E_{\bmxi \sim D(\bx)} \left[\|\nabla_{\bx} f(\bx,\bmxi)\|+|f(\bx,\bmxi)| \left\|\frac{\nabla_{\bx} \mpr(\bmxi\mid\bx)}{\mpr(\bmxi\mid\bx)}\right\| \right]\right\|
    \le \left\|\E_{\bmxi \sim D(\bx)} \left[L_f+f_{\max} M \right]\right\| \\
    &=L_f+f_{\max} M,
    \end{align*}
where the first equality follows from Lemma \ref{sfo_lemma} with $\delta=0$, and the second inequality comes from conditions (i) and (iii) of Assumption \ref{assumption: pvt_v2} and Assumption \ref{asp:compact_set}.
  \end{proof}

\subsection{Proof of Proposition \ref{prop:ogd_regret}}
\begin{proof}
We show that our problem satisfies the assumptions of \cite[Lemma C-5]{besbes2015non}. First, we show that $\{\psi_k\}_{k=1}^R$ in our problem is included in $\mathcal F_s$ defined by \citep[Section 5]{besbes2015non}. $\mathcal F_s$ is a class of sequences $\{g_k\}_{k=1}^R$ of convex cost functions from $\mathcal Y \subset \R^d$ into $\R$, where $\mathcal Y$ is convex, compact, and non-empty. Moreover, $\mathcal F_s$ and $\mathcal Y$ satify the following conditions for all $k \in [R]$:
\begin{itemize}
\item[1.] There is a finite number $G>0$ such that $|g_k(\by)| \le G$ and $\|\nabla_{\by} g_k(\by)\| \le G$ for all $\by \in \mathcal Y$.
\item[2.] There is some $\nu>0$ such that $\{\by \in \R^d: \|\by-\by^*_k\| \le \nu \} \subset \mathcal Y$, where $\by^*_k\in \arg \min_{\by \in \mathcal Y} \psi_k(\by)$.
\item[3.] There are finite numbers $H>0$ and $G>0$ such that $H \bm{I}_d\preceq \nabla_{\by}^2 g_k(\by) \preceq G \bm{I}_d$, where $\bm{I}_d$ is the $d$-dimensional identity matrix.
\end{itemize}
We consider the case of $d=1$ and $\mathcal Y=[-f_{\max}-\kappa, f_{\max}+\kappa]$. Accordingly, $\{\psi_k\}_{k=1}^R$ in our problem is included in $\mathcal F_s$ since the following holds for any $\delta \in \mathcal Y$ and $k \in [R]$:
\begin{align*}
&|\psi_k(\delta)|=\frac{1}{2}(\delta-\E_{\bmxi\sim D(\bx_k)}[f(\bx_k,\bmxi)])^2 \le \frac{1}{2}(2f_{\max}+\kappa)^2,\\
&|\nabla_{\delta} \psi_k(\delta)|=|\delta-\E_{\bmxi\sim D(\bx_k)}[f(\bx_k,\bmxi)]|\le 2f_{\max}+ \kappa,\\
&\arg \min_{\delta'} \psi_k(\delta') \in [-f_{\max},f_{\max}],\ \textrm{and}\\
& \nabla_{\delta}^2 \psi_k(\delta)=1.
\end{align*}

Here, (C-9) in the proof of Lemma C-5 of \cite{besbes2015non} holds by letting $\phi^1(\delta_k,\psi_k):=\delta_{k}-v_k$ since $\E_{v_k}[\delta_{k}-v_k] = \delta_k - \E_{\bmxi\sim D(\bx_k)}[f(\bx_k,\bmxi)] = \nabla \psi_k(\delta_k)$ and $\E_{v_k}[|\delta_{k}-v_k|^2] \le |2f_{\max}+\kappa|^2$. 
Proposition~\ref{prop:ogd_regret} follows from the same argument as in the proof of \cite[Lemma C-5]{besbes2015non}.
\end{proof}

\subsection{Proof of Theorem \ref{thm:convergence}}
\begin{proof}
 When $\delta_{k} \in [-f_{\max},f_{\max}]$, the output $\delta_{k+1}$ of Algorithm \ref{alg:rsag} is included in $[-f_{\max},f_{\max}]$ from $\frac{1}{m_k} { \sum_{\ell=1}^{m_{k}}}f(\bx_{k}^{md},\bmxi^\ell) \in [-f_{\max},f_{\max}]$, $\zeta_{k+1}\in (0,1)$, and line \ref{line:calc_baseline} of Algorithm \ref{alg:rsag}. Therefore, $\delta_k \in [-f_{\max},f_{\max}]$ for all $k \in [R]$ from $\delta_1 \in [-f_{\max},f_{\max}]$.
 From Lemmas \ref{sfo_lemma}-\ref{lem:entire_derivative_bound}, Assumption \ref{asp:compact_set}, and \citep[Corollary 6]{ghadimi2016accelerated}, we have
    \begin{align*}
    \E[\| \mathcal{G}(\bx_R^{md},{\beta_R}) \|^2] \le 96 L_{Ef}\left[ \frac{4L_{Ef} \|\bx_0-\bx^* \|^2}{N(N+1)(N+2)} +\frac{L_{Ef}(\|\bx^*\|^2+H^2)+2\tilde{D}^2}{N}\right].
\end{align*}
Here, in \citep[Corollary 6]{ghadimi2016accelerated}, we let $L_{\psi}:=L_{Ef}$, $L_{f}:=L_{Ef}$, and $\sigma^2:=(L_f+2f_{\max}M)^2$.
Then, to obtain an $\epsilon$-stationary point, we need the iteration number $\hat{N}$ such that 
\begin{align}
    96 L_{Ef}\left[ \frac{4L_{Ef} \|\bx_0-\bx^* \|^2}{\hat{N}(\hat{N}+1)(\hat{N}+2)} +\frac{L_{Ef}(\|\bx^*\|^2+H^2)+2\tilde{D}^2}{\hat{N}}\right] \le \epsilon^2. \label{eq:le_epsilon}
\end{align}
Eq. \eqref{eq:le_epsilon} can be reformulated as
\begin{align*}
    \hat{N}(\hat{N}+1)(\hat{N}+2)\ge \frac{384L_{Ef}^2}{\epsilon^2} \|\bx_0-\bx^* \|^2 + \frac{96 L_{Ef}}{\epsilon^2}(L_{Ef}(\|\bx^*\|^2+H^2)+2\tilde{D}^2)(\hat{N}+1)(\hat{N}+2).
\end{align*}
Therefore, the sufficient condition for \eqref{eq:le_epsilon} is as follows:
\begin{align*}
    \hat{N}^3 \ge \frac{768L_{Ef}^2}{\epsilon^2} \|\bx_0-\bx^* \|^2, \  \hat{N} \ge \frac{192 L_{Ef}}{\epsilon^2}(L_{Ef}(\|\bx^*\|^2+H^2)+2\tilde{D}^2).
\end{align*}
  \end{proof}

\subsection{Proof of Proposition \ref{lem:mpp_ass23}}
\begin{proof}
Assumption \ref{asp:revenue_average} holds since $c(\bmxi)$ is continuous from the definition and \\
$\E_{\bmxi \sim D(\bx)}[s(\bx,\bmxi)] = \E_{\bmxi \sim D(\bx)}\left[\sum_{i=1}^n x_i \xi_i\right]= \sum_{i=1}^n x_i \E_{\bmxi \sim D(\bx)}[\xi_i] 
    = s(\bx,\E_{\bmxi \sim D(\bx)}[\bmxi])$.
Moreover, since $p_k(\bx) \neq 0$ for all $\bx \in \mR^n$ and $k\in\{0,\dots, n\}$ from the definition of $p_i(\bx)$ for each $i\in \{0,1,\dots, n\}$,
%\begin{align*}
$\left(\frac{\nabla_{\bp}\phi(\bp(\bx),\bmxi)}{\phi(\bp(\bx),\bmxi)}\right)_k=\frac{\partial \prod_{i=0}^n {}_mC_{\xi_i} p_i(\bx)^{\xi_i}}{\partial p_k} \frac{1}{\prod_{i=0}^n {}_mC_{\xi_i} p_i(\bx)^{\xi_i}} 
 = \prod_{i=0}^n {}_mC_{\xi_i} p_i(\bx)^{\xi_i}  \frac{\xi_k}{p_k(\bx)}\frac{1}{\prod_{i=0}^n {}_mC_{\xi_i} p_i(\bx)^{\xi_i}} = \frac{\xi_k}{p_k(\bx)}$.  
\end{proof}

\subsection{Proof of Proposition \ref{lem:cphl_asp23}}
\begin{proof}
Assumption \ref{asp:revenue_average} holds since $c(\bmxi)$ is continuous and $s(\bx,\bmxi)=0$.
Moreover, since $0 < p_i(x) < 1$ for all $x \in \mR$ from the definition of $p_i$, we have
\small
\begin{align*}
&\left(\frac{\nabla_{\bp}\phi(\bp(\bx),\bmxi)}{\phi(\bp(\bx),\bmxi)}\right)_k=\frac{\partial \prod_{i \in I} 
{}_{d_i}C_{\xi_i} p_i(x_i)^{\xi_i}(1-p_i(x_i))^{d_i-\xi_i}}{\partial p_k}  \frac{1}{\prod_{i \in I} {}_{d_i}C_{\xi_i} p_i(x_i)^{\xi_i}(1-p_i(x_i))^{d_i-\xi_i}} \\
& = \prod_{i \in I} {}_{d_i}C_{\xi_i}p_i(x_i)^{\xi_i}(1-p_i(x_i))^{d_i-\xi_i} \left(\frac{\xi_k}{p_k(x_k)}-\frac{d_k-\xi_k}{1-p_k(x_k)}\right)\frac{1}{\prod_{i \in I} {}_{d_i}C_{\xi_i} p_i(x_i)^{\xi_i}(1-p_i(x_i))^{d_i-\xi_i}}\\
& = \frac{\xi_k}{p_k(x_k)}-\frac{d_k-\xi_k}{1-p_k(x_k)}. 
\end{align*}
\normalsize
\end{proof}

\subsection{Proof of Lemma \ref{lem:sfo_accel}}

\begin{proof}
It follows from the definition of $\phi$ that %\takeda{The definition of $\phi$ is in Proposition 12. The same $\phi$ used in Algorithm 1??}\memo{Algorithm 1の方の記号を変更しました．}
\begin{align}
 &\E_{\bmxi \sim D(\bx)}\left[\delta\frac{\nabla_{\bx}\phi(\bp(\bx),\bmxi)}{\phi(\bp(\bx),\bmxi)} \right] = \int_{\bmxi \in \Xi} \left( \delta \frac{\nabla_{\bx}\phi(\bp(\bx),\bmxi)}{\phi(\bp(\bx),\bmxi)}\right) d\mpr(\bmxi\mid\bx) =\delta \int_{\bmxi \in \Xi} \nabla_{\bx} \mpr(\bmxi\mid\bx) d\bmxi \nonumber \\
  &=\delta \nabla_{\bx} \int_{\bmxi \in \Xi} \mpr(\bmxi\mid\bx) d\bmxi=\delta \nabla_{\bx}(1)=0, \label{eq:cost_case_delta_int_0}
\end{align}
where the third equality comes from Lemma \ref{lem:int_lim_change_2} with $q(\bmxi) = 1$.
Since 
\begin{align*}
\E_{\bmxi \sim D(\bx)}[f(\bx,\bmxi)] = -s(\bx,\E_{\bmxi \sim D(\bx)}[\bmxi]) + \int_{\bmxi \in \Xi} c(\bmxi) \phi(\bp(\bx),\bmxi) d\bmxi
\end{align*}
from Assumptions \ref{asp:revenue_average} and~\ref{asp:choice_model}, we have
 \begin{align*}
 \nabla_{\bx} \E_{\bmxi \sim D(\bx)}[f(\bx,\bmxi)] 
 &= \nabla_{\bx} \left( -s(\bx,\E_{\bmxi \sim D(\bx)}[\bmxi]) + \int_{\bmxi \in \Xi} c(\bmxi) \phi(\bp(\bx),\bmxi) d\bmxi\right)\\ 
 & = - \nabla_{\bx} s(\bx,\E_{\bmxi \sim D(\bx)}[\bmxi]) + \nabla_{\bx} \int_{\bmxi \in \Xi} c(\bmxi) \phi(\bp(\bx),\bmxi) d\bmxi\\
 &= -\nabla_{\bx} s(\bx,\E_{\bmxi \sim D(\bx)}[\bmxi])  +  \int_{\bmxi \in \Xi} \nabla_{\bx} \left( c(\bmxi) \phi(\bp(\bx),\bmxi) \right) d\bmxi \\
 &= -\nabla_{\bx} s(\bx,\E_{\bmxi \sim D(\bx)}[\bmxi]) + \E_{\bmxi \sim D(\bx)}\left[c(\bmxi)\frac{\nabla_{\bx}\phi(\bp(\bx),\bmxi)}{\phi(\bp(\bx),\bmxi)} \right] \\
 &= -\nabla_{\bx} s(\bx,\E_{\bmxi \sim D(\bx)}[\bmxi]) + \E_{\bmxi \sim D(\bx)}\left[(c(\bmxi)-\delta)\frac{\nabla_{\bx}\phi(\bp(\bx),\bmxi)}{\phi(\bp(\bx),\bmxi)} \right] \\
&= -\nabla_{\bx} s(\bx,\E_{\bmxi \sim D(\bx)}[\bmxi])  + \E_{\bmxi \sim D(\bx)}\left[(c(\bmxi)-\delta) \frac{d\bp(\bx)}{d\bx}\frac{\nabla_{\bp}\phi(\bp(\bx),\bmxi)}{\phi(\bp(\bx),\bmxi)} \right].
\end{align*}
Here, the third equality holds when $\bmxi$ is a discrete random vector.
When $\bmxi$ is a continuous random vector, the third equality comes from Assumption \ref{asp:continuous_dist} and Lemma \ref{lem:int_lim_change_2} by letting $q(\bmxi):=c(\bmxi)$.
The fourth equality  is due to the fact that $\phi(\bp(\bx),\bmxi)=\mpr(\bmxi\mid \bx)\neq 0$ from condition (ii) of Assumption \ref{assumption: pvt_v2}.
The fifth equality comes from \eqref{eq:cost_case_delta_int_0}.
  \end{proof}

\subsection{Proof of Lemma \ref{lemma:SFO_variance_2}}
\begin{proof}
Under Assumptions \ref{asp:revenue_average} and \ref{asp:choice_model}, we have 
\begin{align} 
\nabla_{\bx} \E_{\bmxi \sim D(\bx)}[f(\bx,\bmxi)] = -\nabla_{\bx}s(\bx,\E_{\bmxi \sim D(\bx)}[\bmxi]) + \nabla_{\bx}\int_{\bmxi \in \Xi} c(\bmxi) \phi(\bp(\bx),\bmxi) d\bmxi. \label{eq:grad_obj}
\end{align}
From Lemma \ref{lem:sfo_accel},
 \begin{align}
 \nabla_{\bx} \E_{\bmxi \sim D(\bx)}[f(\bx,\bmxi)] 
= -\nabla_{\bx} s(\bx,\E_{\bmxi \sim D(\bx)}[\bmxi])  + \E_{\bmxi \sim D(\bx)}\left[(c(\bmxi)-\delta) \frac{d\bp(\bx)}{d\bx}\frac{\nabla_{\bp}\phi(\bp(\bx),\bmxi)}{\phi(\bp(\bx),\bmxi)} \right].\label{eq:grad_g3}
\end{align}
Then, from \eqref{eq:grad_obj} and \eqref{eq:grad_g3},
\begin{align}
\E_{\bmxi \sim D(\bx)}\left[(c(\bmxi)-\delta)\frac{d\bp(\bx)}{d\bx} \frac{\nabla_{\bp} \phi(\bp(\bx),\bmxi)}{\phi(\bp(\bx),\bmxi)}\right] = \nabla_{\bx} \int_{\bmxi \in \Xi} c(\bmxi) \phi(\bp(\bx),\bmxi) d\bmxi. \label{eq:letter}
\end{align}
Then,
\begin{align*}
&\E_{\bmxi' \sim D(\bx)}[\| g_2(\bx,\bmxi') - \nabla_{\bx} \E_{\bmxi \sim D(\bx)}[f(\bx,\bmxi)]\|^2] \\
&=  \E_{\bmxi' \sim D(\bx)}\left[\left\| -\nabla_{\bx} s(\bx,\E_{\bmxi \sim D(\bx)}[\bmxi])  + (c(\bmxi')-\delta) \frac{d\bp(\bx)}{d\bx} \frac{\nabla_{\bp} \phi(\bp(\bx),\bmxi')}{\phi(\bp(\bx),\bmxi')} \right. \right. \\
& \qquad \qquad \qquad \left. \left. + \nabla_{\bx} s(\bx,\E_{\bmxi \sim D(\bx)}[\bmxi]) - \nabla_{\bx} \int_{\bmxi \in \Xi} c(\bmxi) \phi(\bp(\bx),\bmxi) d\bmxi\}\right\|^2\right] \\
&= \E_{\bmxi' \sim D(\bx)}\left[\left\| (c(\bmxi')-\delta) \frac{d\bp(\bx)}{d\bx} \frac{\nabla_{\bp} \phi(\bp(\bx),\bmxi')}{\phi(\bp(\bx),\bmxi')} - \nabla_{\bx} \int_{\bmxi \in \Xi} c(\bmxi) \phi(\bp(\bx),\bmxi) d\bmxi\}\right\|^2\right]\\
&= \E_{\bmxi' \sim D(\bx)}\left[\left\| (c(\bmxi')-\delta) \frac{d\bp(\bx)}{d\bx} \frac{\nabla_{\bp} \phi(\bp(\bx),\bmxi')}{\phi(\bp(\bx),\bmxi')}\right\|^2\right] - \left\|\nabla_{\bx} \int_{\bmxi \in \Xi} c(\bmxi) \phi(\bp(\bx),\bmxi) d\bmxi\}\right\|^2\\
&\le \E_{\bmxi' \sim D(\bx)}\left[\left\|(c(\bmxi')-\delta) \frac{d\bp(\bx)}{d\bx} \frac{\nabla_{\bp} \phi(\bp(\bx),\bmxi')}{\phi(\bp(\bx),\bmxi')}\right\|^2\right]=\E_{\bmxi' \sim D(\bx)}\left[\left\| (c(\bmxi')-\delta) \frac{\nabla_{\bx} \mpr(\bmxi' \mid \bx)}{\mpr(\bmxi' \mid \bx)}\right\|^2\right] \\
&\le \E_{\bmxi' \sim D(\bx)}\left[(2c_{\max})^2M^2 \right] 
= 4(c_{\max}M)^2,
\end{align*}
where the first equality comes from \eqref{eq:grad_obj}, and the third equality follows from \eqref{eq:letter}.
The second inequality follows from condition (iii) of Assumption \ref{assumption: pvt_v2}.
\end{proof}

\subsection{Proof of Theorem \ref{thm:convergence_multagent}}
\begin{proof}
When $\delta_{k} \in [-c_{\max},c_{\max}]$, the output $\delta_{k+1}$ of Algorithm \ref{alg:rsag_specialized} is included in $[-c_{\max},c_{\max}]$ from $ \frac{1}{m_k} \sum_{\ell=1}^{m_{k}}c(\bmxi^\ell) \in [-c_{\max},c_{\max}]$, $\zeta_{k+1}\in (0,1)$, and the update rule for $\delta_k$ in Algorithm \ref{alg:rsag_specialized}. 
Therefore, $\delta_k \in [-c_{\max},c_{\max}]$ for all $k \in [R]$ from $\delta_1 \in [-c_{\max},c_{\max}]$.
From Lemmas \ref{lem:entire_derivative_bound}, \ref{lem:sfo_accel}, \ref{lemma:SFO_variance_2}, and  \citep[Corollary 6]{ghadimi2016accelerated}, we have 
\begin{align*}
\E[\| \mathcal{G}(\bx_R^{md},{\beta_R}) \|^2] \le 96 L_{Ef}\left[ \frac{4L_{Ef} \|\bx_0-\bx^* \|^2}{N(N+1)(N+2)} +\frac{L_{Ef}(\|\bx^*\|^2+H^2)+2\tilde{D}^2}{N}\right].
\end{align*}
Here, in \citep[Corollary 6]{ghadimi2016accelerated}, we let $L_{\psi}:=L_{Ef}$, $L_{f}:=L_{Ef}$, and $\sigma^2:=4(c_{\max}M)^2$.
Then, as in the proof of Theorem \ref{thm:convergence}, we need the iteration number $\hat{N}$ to obtain an $\epsilon$-stationary point such that:
\begin{align*}
    \hat{N}^3 \ge \frac{768L_{Ef}^2}{\epsilon^2} \|\bx_0-\bx^* \|^2, \  \hat{N} \ge \frac{192 L_{Ef}}{\epsilon^2}(L_{Ef}(\|\bx^*\|^2+H^2)+2\tilde{D}^2).
\end{align*}
  \end{proof}

\section{Details of our experiments} \label{app:detail_exp}
\subsection{Common Settings}
All experiments were conducted on a computer with an AMD EPYC 7413 24-Core Processor, 503.6 GiB of memory RAM, and Ubuntu 20.04.6 LTS. The program code was implemented in Python 3.8.3.

\subsection{Settings of Baselines}
\textbf{L2-Regularized Repeated Gradient Descent (L2-RGD($\alpha$)):} This method is described in Section 2.2. We used the fixed step size $\eta_k:=0.01$ at each iteration $k$.
\\
\textbf{Bayesian Optimization (BO):} We used GPyOpt, a Python open-source library for Bayesian optimization \citep{gpyopt2016}. We used the default setting of the library for parameters other than the termination criteria.
\\
\textbf{Simultaneous Perturbation Stochastic Approximation (SPSA):} At each iteration $k$, this method updates the current iterate by using the stochastic perturbation gradient:
\begin{align*}
\frac{f(\bx^k+c_k\Delta^k,\bmxi^{k,1})-f(\bx^k-c_k\Delta^k,\bmxi^{k,2})}{c_k\Delta^k},
\end{align*}
where $c_k:=\frac{1}{(k+1)^{0.101}}$, each element of $\Delta^k$ is sampled from a Rademacher distribution (i.e. Bernoulli $\pm 1$ with probability $0.5$), and $\bmxi^{k,1}$ and $\bmxi^{k,2}$ are random vectors sampled from the distribution $D(\bx^k)$. We set $a_k:=\frac{0.16}{(100+k+1)^{0.602}}$ as the stepsize at each iteration. The settings of $c_k$, $\Delta^k$, and $a_k$ are based on \citep[Section III]{spall1998implementation}.
\\ 
\textbf{ Projected Sub-gradient Descent for Average Demand (PSD-AD):} This method is a projected subgradient descent method for 
\begin{align*}
\min_{\bx \in [x_{\min},x_{\max}]^n} (-s(\bx,\bar{\bmxi}(\bx)) + c(\bar{\bmxi}(\bx))),
\end{align*}
where $\bar{\bmxi}(\bx):=\E_{\bmxi \sim D(\bx)}[\bmxi]$, which represents the average demand for $\bx$. We set the step size at each iteration so that the objective value decreases by repeatedly multiplying by $\delta=0.9$.

\end{document}